\documentclass[12pt]{article}
\usepackage[small,compact]{titlesec}
\usepackage{amsmath,amssymb,amsfonts,setspace} 
\usepackage{graphicx,caption,epsfig,subfigure,epsfig,wrapfig}
\usepackage{url,color,verbatim,algorithmic}
\usepackage[sort,compress]{cite}
\usepackage[ruled,boxed]{algorithm}
\usepackage[scaled]{helvet}
\usepackage[T1]{fontenc}
\usepackage{booktabs,multirow} 
\usepackage[bookmarks=true, bookmarksnumbered=true, colorlinks=true,   pdfstartview=FitV,
linkcolor=blue, citecolor=blue, urlcolor=blue]{hyperref}

\let\OLDthebibliography\thebibliography
\renewcommand\thebibliography[1]{
  \OLDthebibliography{#1}
  \setlength{\parskip}{0pt}
  \setlength{\itemsep}{0pt}
}
\usepackage{authblk} 

\def\wv{\widehat{v}} 
\def\hatu{\widehat{u}}  

\def\bA{\mathbf{A}}    
\def\bb{\mathbf{b}}

\def\R{\mathbb{R}}       \def\C{\mathbb{C}}

\newtheorem{theorem}{Theorem}

\newtheorem{remark}[theorem]{Remark}

\numberwithin{equation}{section}
\numberwithin{theorem}{section}

\usepackage{soul}

\definecolor{greenrb}{rgb}{0.2,0.6,0.2}

\title{Shock trace prediction by reduced models for \\ a viscous stochastic Burgers equation}

\author[1]{Nan Chen}
\author[2]{Honghu Liu}
\author[3]{Fei Lu}
\affil[1]{\footnotesize Department of Mathematics, University of Wisconsin-Madison, Madison, WI 53705, USA (\href{chennan@math.wisc.edu}{chennan@math.wisc.edu})}
\affil[2]{\footnotesize  Department of Mathematics, Virginia Tech, Blacksburg, VA 24061, USA (\href{hhliu@vt.edu}{hhliu@vt.edu})}
\affil[3]{\footnotesize  Department of Mathematics, Johns Hopkins University, Baltimore, MD 21218, USA (\href{feilu@math.jhu.edu}{feilu@math.jhu.edu})}

\date{\today}

\begin{document}\maketitle

\begin{abstract}
Viscous shocks are a particular type of extreme events in nonlinear multiscale systems, and their representation requires small scales.
Model reduction can thus play an important role in reducing the computational cost for an efficient prediction of shocks. Yet, reduced models typically aim to approximate large-scale dominating dynamics, which do not resolve the small scales by design. To resolve this representation barrier, we introduce a new qualitative characterization of the space-time locations of shocks, named as the ``shock trace'', via a space-time indicator function based on an empirical resolution-adaptive threshold. Different from the exact shocks, the shock traces can be captured within the representation capacity of the large scales, which facilitates the forecast of the timing and locations of the shocks utilizing reduced models. Within the context of a viscous stochastic Burgers equation, we show that a data-driven reduced model, in the form of nonlinear autoregression (NAR) time series models, can accurately predict the random shock traces, with relatively low rates of false predictions. The NAR model significantly outperforms the corresponding Galerkin truncated model in the scenario of either noiseless or noisy observations. The results illustrate the importance of the data-driven closure terms in the NAR model, which account for the effects of the unresolved small scale dynamics on the resolved ones due to nonlinear interactions.  
\end{abstract}


\tableofcontents

\section{Introduction}

Extreme events occur in many high-dimensional multiscale systems in geophysics, engineering, neural science and material science \cite{ghil2011extreme, kwasniok2012data, sheard2009principles, saha2017extreme,wilcox1988multiscale, majda2016introduction, tao2009multiscale, majda2018model}. Prediction of them with uncertainty quantification has significant scientific and societal impacts. However, it can be computationally prohibitive to run these systems for the ensemble prediction that aims to quantify the uncertainty.  Reduced models can bring down the computational cost by orders of magnitude compared with that for the original system, 
providing surrogate models that make the ensemble prediction feasible. 
In particular, stochastic reduced models have been built to quantify the uncertainty 
 \cite{palmer2010stochastic, majda2012lessons, smith2013uncertainty,CL15,LTC17} and have been applied to  
 forecast extreme events \cite{kwasniok2012data,majda2018strategies, cousins2014quantification, qi2018predicting, franzke2012predictability}. Many of these reduced models primarily aim at the statistical forecast of the extreme values, for example, predicting the probability density functions.  Other reduced models are designed for the short-term prediction of intermittent time series from nature that contain extreme events \cite{chen2020predicting, cousins2016reduced}. In addition, predicting the onset of extreme events in complex systems via suitable statistical or machine learning tools have also been developed and applied to several important turbulent systems \cite{qi2020using, sapsis2021statistics, farazmand2019extreme,wan2018data}.

Shock (or shock wave) is a special type of extreme events that is observed in many complex nonlinear systems \cite{courant1999supersonic, shapiro1953dynamics, anderson2010fundamentals}. Shocks have a unique feature that makes them extremely difficult to predict: a shock is characterized by an abrupt change of the state, with the spatial derivative rather than the state itself reaching an extreme value. Thus, the exact representation of a shock requires the small fast scales (see Figure \ref{fig:pdf_Du} and discussions in Section \ref{sec:shockTrace}), which are beyond the reach of a reduced model that models the large scales. Yet, the large scales contain rich information about the shocks, because they dominate the small scales through the nonlinear interactions. Thus, it is of interest to uncover the relation between the shocks and the large scales, and to investigate the possibility of predicting partial information of shocks by the reduced models. As far as we know,  despite the importance of predicting the shocks, this issue has not been addressed yet.

This paper shows that reduced models can predict the timing and location of random shocks of viscous stochastic Burgers equations \cite{e1997_ProbabilityDistributiona, e2000_InvariantMeasures}. The Burgers equation is a prototype for nonlinear conservation equations that can develop shocks. Here, on top of the nonlinear deterministic dynamics, external stochastic forces are added into the Burgers equation. These stochastic forces are smooth in space and white in time. They play an important role in randomly triggering shocks. An additional small viscous term is further incorporated into the equation such that the random shocks will be dissipated before the appearance of a discontinuity. Therefore, these shocks are known as viscous shocks \cite{dunlap2021viscous}. In addition to the fundamental mathematical topics such as the existence of solution and the invariant measure \cite{DaPrato1994stochastic,e1997_ProbabilityDistributiona,e2000_InvariantMeasures}, the statistics of shocks has been studied for  the Burgers equation in various contexts, either with viscosity or in the inviscid limit, with or without stochastic forcing, and possibly subject to additional random initial conditions \cite{sinai1992statisticsShocks,dunlap2021viscous,cho2014_StatisticalAnalysis,bec2007_BurgersTurbulence,chorin2005_ViscositydependentInertial}. Burgers equation has also frequently served as a testbed for illustrating various model reduction techniques, although usually not in the context of shock tracking; see {\it e.g.}~\cite{dolaptchiev2013stochastic,wangTwolevelDiscretizations2011,mou2020dd,beck2009_UsingGlobal}.

To effectively predict the spatiotemporal structures of the shocks, we introduce a new qualitative characterization via a simple indicator function that depends on both time and space. It is named as ``shock trace'' (see Section \ref{sec:shockTrace}). By prescribing an empirical threshold value, the spatial derivative of the reconstructed spatiotemporal solution from a reduced model is mapped to this indicator function. If the spatial derivative of the forecast solution is more negative than the threshold, then it is mapped to 1, indicating occurrence of a shock. Otherwise, it is mapped to 0, standing for no shock occurrence at that location. This new characterization has several desirable features in facilitating the prediction of the shocks. First, it mitigates the complexity of the shock forecast as it avoids the forecast of the precise values of the solution. Nevertheless, the simplified representation of the solution via the indicator function preserves the entire spatiotemporal structure of the shocks, providing an effectual indicator of the abrupt and sharp changes in the solution. Second, the threshold value of the model solution that maps to the indicator function is adaptive to spatial resolutions (i.e., the number of Fourier modes being adopted). Therefore, the same mapping criterion can be applicable to the solutions of both the full model and reduced models. Third, the indicator function provides a simple but effective way to quantitatively count the false predictions, allowing a quantification of the forecast uncertainty.

We consider a data-driven reduced model in the form of a nonlinear autoregression (NAR) model \cite{Lu20Burgers,LinLu21}. The NAR model has three attractive features as a prototype reduced model. First, it includes only the time evolution of the large-scale Fourier coefficients with low wave-numbers, and it allows a large time-step size since there is no stiffness in these large scales. Therefore, it significantly decreases the computational cost, often several orders of magnitude from an accurate full model. Second, as a closure model, the NAR model effectively parameterizes the nonlinear feedback from the unresolved small scales, and its parametric form is derived from a numerical integrator of the system. Third, its parameters can be efficiently estimated by least squares regression from data (see Section \ref{sec:NAR}).   

Numerical results show that the NAR model can predict the shock traces accurately with low rates of false predictions, almost as good as the corresponding projection of the full model solution. The NAR model significantly outperforms the Galerkin truncated system, either in predictions from true initial conditions or in data assimilation from noisy observations, highlighting the importance and effectiveness of the data-driven closure that captures the effects of unresolved scales; see Section \ref{sec:numerics}.

The rest of the paper is organized as follows. Following a brief review of the basic properties of the stochastic Burgers equation and its numerical integration scheme, we introduce the concept of shock trace based on the indicator function with a resolution-adaptive threshold in Section \ref{sec:SBE}. Section \ref{sec:NAR} presents the NAR modeling framework and its parameter inference. The forecast skill of the reduced model is studied in Section \ref{sec:numerics} for cases with both weak and strong stochastic forces. Conclusion and some final remarks are then presented in Section \ref{sec:conclusion}.

\section{The Viscous Stochastic Burgers Equation and Shock Trace} \label{sec:SBE} 
\subsection{The viscous stochastic Burgers equation}
The model considered in this article is the following viscous stochastic Burgers equation posed on $(0, 2\pi)$ supplemented with periodic boundary conditions and suitable initial condition \cite{e1997_ProbabilityDistributiona, e2000_InvariantMeasures}:
\begin{eqnarray} \label{SBE}
\begin{aligned}
&\partial_{t} u = \nu \partial_{xx} u-u \partial_{x} u+ f(x,t), \quad 0 < x <2 \pi, \; t>0, \\
 &u(0,t)=u(2\pi ,t) , \quad  \partial_{x} u(0,t)= \partial_{x} u(2\pi ,t), \quad  t \ge 0, \\
 & u(x,0) = u_0(x), \quad 0 <  x  < 2 \pi.
\end{aligned}
\end{eqnarray}
Here $\nu>0$ is the viscosity constant, $\sigma>0$
represents the strength of the stochastic force, and $u_0$ is a given square-integrable function. The stochastic force $f(x,t)$ is smooth in space and white in time \cite{gardiner2004handbook}, acting on a few low-frequency (i.e., large-scale) Fourier modes, which is chosen to be of the form
\begin{equation} \label{eq:sforce}
f(x,t)  =\sigma \sum_{m=1}^{K_0}  \Big( \sin(mx) \dot{W}_m(t)+ \cos(mx) \dot{W'}_m(t) \Big),
\end{equation}
where $\{W_m,W'_m\}$ are independent Brown motions with $\{\dot{W}_m,\dot{W}'_m\}$ denoting the white noises and $K_0$ is a fixed positive integer. With the chosen boundary conditions and the above form of the stochastic force $f$, one can check that the quantity $\int_0^{2\pi } u(x,t) \, \mathrm{d}x$ is conserved for all $t \ge 0$. Without loss of generality, we assume that the initial condition has mean zero, leading thus to
\begin{equation} \label{Eq_meanzero_cond}
\int_0^{2\pi } u(x,t) \, \mathrm{d}x=0, \quad t \ge 0.
\end{equation}

The equation \eqref{SBE} is interpreted as an infinite dimensional stochastic differential equation (SDE) for the corresponding Fourier modes,
\begin{eqnarray}
\frac{\mathrm{d}}{\mathrm{d} t}\hatu_{k} 
&=&- \nu k^{2}\hatu_{k}- \frac{ik}{2}\sum_{l=-\infty}^{\infty}\hatu_{l}\hatu_{k-l}  + \widehat{f}_k(t) \label{FM}
\end{eqnarray}
where $\hatu_{k}$ are the Fourier coefficients: 
\begin{equation*}
\hatu_{k}(t)=\mathcal{F}[u]_{k}=\frac{1}{2\pi }%
\int_{0}^{2\pi }u(x,t)e^{-ikx}dx,\,\,\,u(x,t)=\mathcal{F}^{-1}[\hatu%
]=\sum_{k=-\infty}^\infty \hatu_{k}(t)e^{ikx},
\end{equation*}%
with $\mathcal{F}[u]$ being the Fourier transform of $u$. Here $\widehat{f}_k(t)$ are white noises consisting of linear combinations of $\{\dot{W}_k,\dot{W}'_k\}$ in \eqref{eq:sforce} for $|k|\leq K_0$, and $\widehat{f}_k(t) =0$ for $k>K_0$. Note also that $\hatu_0(t) \equiv 0$ thanks to \eqref{Eq_meanzero_cond}.

The above viscous stochastic Burgers equation driven by the stochastic force, which is smooth in space and white in time, has been proven to have an invariant measure in \cite{e1997_ProbabilityDistributiona,e2000_InvariantMeasures}. We refer the readers to \cite{DaPrato1994stochastic} for the Burgers equation driven by spatiotemporal white noises. Note that the ``one force, one solution'' principle holds: for each realization of the stochastic force, there exists a unique solution globally and the random attractor consists of a single trajectory, almost surely. For each realization of the force, the shocks formulate randomly, but dissipate before reaching a discontinuity due to the viscosity. Thus, they are called viscous shocks \cite{dunlap2021viscous}.

\subsection{Numerical integration scheme of the full system and parameter regimes}\label{sec:num}
A high-resolution numerical integration scheme is adopted for the full system \eqref{SBE} to generate the true solution.
The numerical scheme utilized here is a Galerkin spectral method \cite{guo1998spectral}. More specifically, the function $u(x,t)$ is represented at grid points $x_{j}=j \Delta x$ with $j=0,\dots ,2N-1$ and $\Delta x = \frac{2\pi }{2N}$. The Fourier transform $\mathcal{F}$ is replaced by discrete Fourier transform
\begin{equation}
\begin{aligned}
\hatu_{k}(t) &=\mathcal{F}_{2N}[u]_{k}=\sum_{j=0}^{2N-1}u(x_{j},t)e^{-ikx_{j}}, \quad k = -N+1, -N+2, \ldots, N,\\
u(x_{j},t) &=\mathcal{F}_{2N}^{-1}[\hatu]_{j}=\frac{1}{2N}%
\sum_{k=-N+1}^{N}\hatu_{k}e^{ikx_{j}}, \quad  j = 0, 1, \ldots ,2N-1.
\end{aligned}
\end{equation}
Since $u$ is real, we have $\hatu_{-k}=\hatu_{k}^{\ast }$. Recall also that $\hatu_{0}$ is identically zero thanks to the mean-zero assumption on the initial data (cf.~\eqref{Eq_meanzero_cond}). Then, by setting $\hatu_{-N}=0$ to simply the notations,  
we obtain a truncated system
\begin{align} \label{DFM}
\frac{\mathrm{d}}{\mathrm{d} t}\hatu_{k}=-\nu k^{2}\hatu_{k}-\frac{ik}{2} \sum_{\substack{|k-l|\leq N, |l|\leq N }}  \hatu_l \hatu_{k-l}+ \widehat{f}_k, \text{ with } |k|=1,\dots ,N.  
\end{align}
The system \eqref{DFM} is solved using the exponential time differencing fourth order Rouge--Kutta method (ETDRK4) (see \cite{CM02, KT05}) with the standard $3/2$ zero-padding for dealiasing (see e.g., \cite{GO77}), where the force term $\widehat{f}_k$ is treated as a constant in each time step. Such a hybrid scheme is of strong order 1, but it has an advantage of preserving both the numerical stability of ETDRK4 and the simplicity of the Euler--Maruyama scheme.

Hereafter, for a given $K < N$, we call $u_K(x,t)= \sum_{|j|\leq K} \hatu_j(t) e^{-ijx}$ the $K$-mode projection of the full model solution $u_N(x,t )= \sum_{|j|\leq N} \hatu_j(t) e^{-ijx}$, where $ \hatu_j(t)$ is solved from the above full model \eqref{DFM} with $N$-pairs of Fourier modes. The reduced model to be introduced later in Section~\ref{sec:NAR} aims to approximate this $K$-mode projection for a suitable $K$.

The parameters used in the numerical experiments are chosen as follows. The viscosity is set to be $\nu=0.02$, which is small enough to allow shocks to emerge, when subject to the stochastic force further specified below. Of course, a smaller viscosity constant usually demands a higher spatiotemporal resolution in order to obtain numerically accurate solutions of \eqref{SBE} and hence an accurate description of the emergence of the shocks. For the chosen $\nu$, we have set $N=128$ for the high-dimensional Galerkin system \eqref{DFM} and we used $\Delta t=0.001$ as the time-step size. Only the first four pairs of Fourier modes are forced with $K_0 = 4$ in \eqref{eq:sforce}; and we consider two values for the strength of the stochastic force, $\sigma=0.2$ and $\sigma=1$, which correspond to two dynamical regimes exhibiting moderate and strong ``turbulent'' behavior, respectively. For the regime with larger $\sigma$, shocks appear more frequently both in time and space, and the (spatial) gradients presented in the viscous shock profiles become sharper as well; cf.~Figures~\ref{fig:pdf_Du} and \ref{fig:traceShock_true}. It is also worth mentioning that, for both forcing scenarios, we checked that the mean Courant--Friedrichs--Lewy (CFL) number (computed over a sufficiently long solution trajectory) is well below one:  it is 0.045 for the case $\sigma = 0.2$ and 0.139 for the case $\sigma = 1$. For latter reference, we summarize the parameters in Table~\ref{tab:RMsettings}, in which we also listed the dimension and time-step size used for the reduced models to be presented in Section~\ref{sec:NAR} below.

\bigskip
 \begin{table}[tbh!]
  \caption{Parameter settings of the full and reduced models}  \label{tab:RMsettings}
  \centering
\begin{tabular}{ c ll }
\toprule
\multirow{5}{*}{Full model} &  $\nu =0.02$ & viscosity constant   \\ 
                                            & $N=128$ & number of modes  \\
                                            & $\Delta t=0.001$ & time-step size\\
                                           & $K_0=4$ & number of modes in the stochastic force\\
                                            & $\sigma = 0.2 \text{ or } 1$ & strength of stochastic force \\[2.5px] 
\hline
\multirow{2}{*}{Reduced  models} & $K=8$     & number of modes in the reduced model \rule{0pt}{2.5ex}\\
                                                    & $\delta  = \Delta t \times 10$  &   time-step size of reduced model \\ 
\bottomrule
\end{tabular}
\end{table}

\subsection{Trace of the random viscous shocks} \label{sec:shockTrace}

A new simple criterion is developed here to qualitatively characterize the occurrence of the random shocks. This criterion is adaptive to spatial resolutions, and is thus applicable to reduced models as well. The procedure of applying this new criterion to trace the shocks is as follows:
\begin{enumerate}
  \item Prescribe a statistical threshold value for the spatial derivative of the process to determine the occurrence of shocks: a shock occurs at $(x,t)$ if the spatial derivative is more negative than the threshold value. The threshold value is computed empirically from a large ensemble of sample trajectories of the process. 
  \item Define an indicator function that identifies the shock trace: it has value 1 at $(x,t)$ if a shock occurs according to the criterion in step 1 above and is 0 otherwise.  
\end{enumerate}
This new criterion mitigates the complexity of the shock forecast as it avoids the forecast of the precise values of the solution.  Nevertheless, the simplified representation of the solution via the indicator function preserves the entire spatiotemporal structure of the shocks, providing an effectual indicator of the abrupt and sharp changes in
the solution. Note that the threshold value is adaptive to spatial resolutions (i.e., the number of Fourier modes being adopted). Therefore, the same mapping criterion is applicable to the solutions of both the full model and reduced models. In addition, the indicator function provides a simple but effective way to quantitatively count the false predictions, allowing a quantification of the forecast uncertainty, as will be seen in Section~\ref{sec:numerics}.

We set the threshold value to be one standard deviation above the mean of the most negative spatial derivatives, which are empirically computed from a large ensemble of trajectories. This criterion applies to processes with different spatial resolutions, e.g., the $k$-mode projection $u_k(x,t) = \sum_{|j|\leq k} \hatu_j(t) e^{ijx} $ of the full model solution, for any $k\geq 1$.   More specifically, 
we compute the empirical threshold value $\tau$ as follows, using this $k$-mode projection $u_k$ as an example. First,  generate a large ensemble of $M$ representative trajectories of the full model whose initial conditions are sampled from a long trajectory of the full model. Denote the $k$-mode projection of these trajectories by $\{u_k^{(m)}(x,t_l), 0\leq l \leq L]\}_{m=1}^M$, where $t_l= l \delta $ and $L=\frac{T}{\delta}$ with $\delta$ corresponding to the sampling frequency and $T$ the length of the time interval over which the solution trajectories are computed. 
Second, compute the threshold value $\tau_k$ for the $k$-mode projection as
\begin{equation}\label{eq:tau}
\tau_k =  \overline{D}_k + \eta_{k}, \quad \text{ with }   \overline{D}_k = \frac{1}{LM}\sum_{l=1}^L \sum_{m=1}^M   D_k^{(m)}(t_l), \,     \eta_k^2 = \frac{1}{LM}\sum_{l=1}^L \sum_{m=1}^M   | D_k^{(m)}(t_l)- \overline{D}_k |^2,
\end{equation}
where $D_k^{(m)} (t_l)= \min_x \partial_{x} u_k^{(m)}(x,t_l)$ denotes the most negative spatial derivative at $t = t_l$ for the $m$-th trajectory of $u_k$. Note that $\overline{D}_k$ and $\eta_k$ are respectively the mean and standard deviation of $D_k^{(m)} (t_l)$ among all times and all trajectories.   

\begin{figure}[htb]
     \centering
\hspace{1mm}		\subfigure[Typical shocks  $\sigma=0.2$]{\includegraphics[scale=.47]{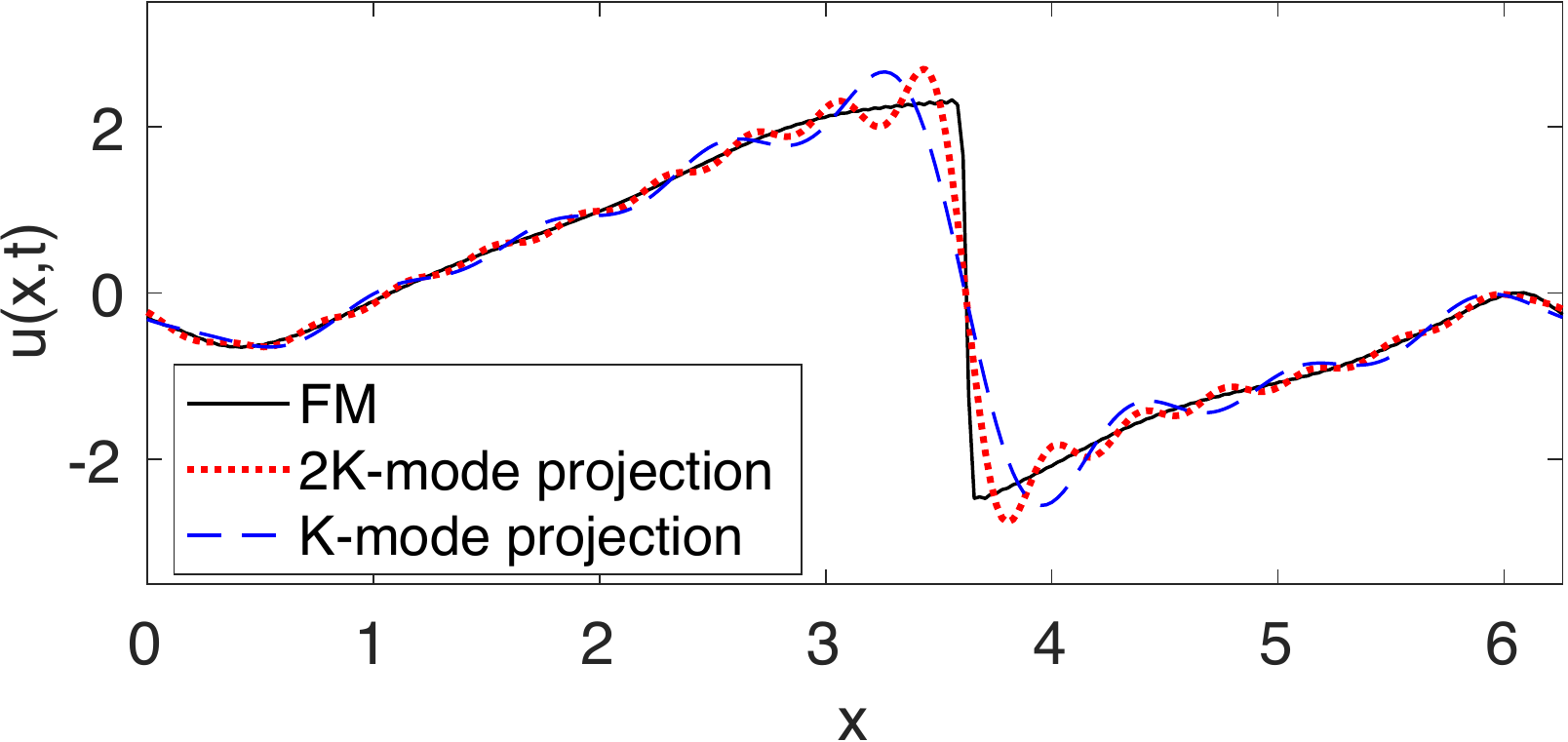}} \hspace{2mm}
		\subfigure[Typical shocks $\sigma=1$]{\includegraphics[scale=.47]{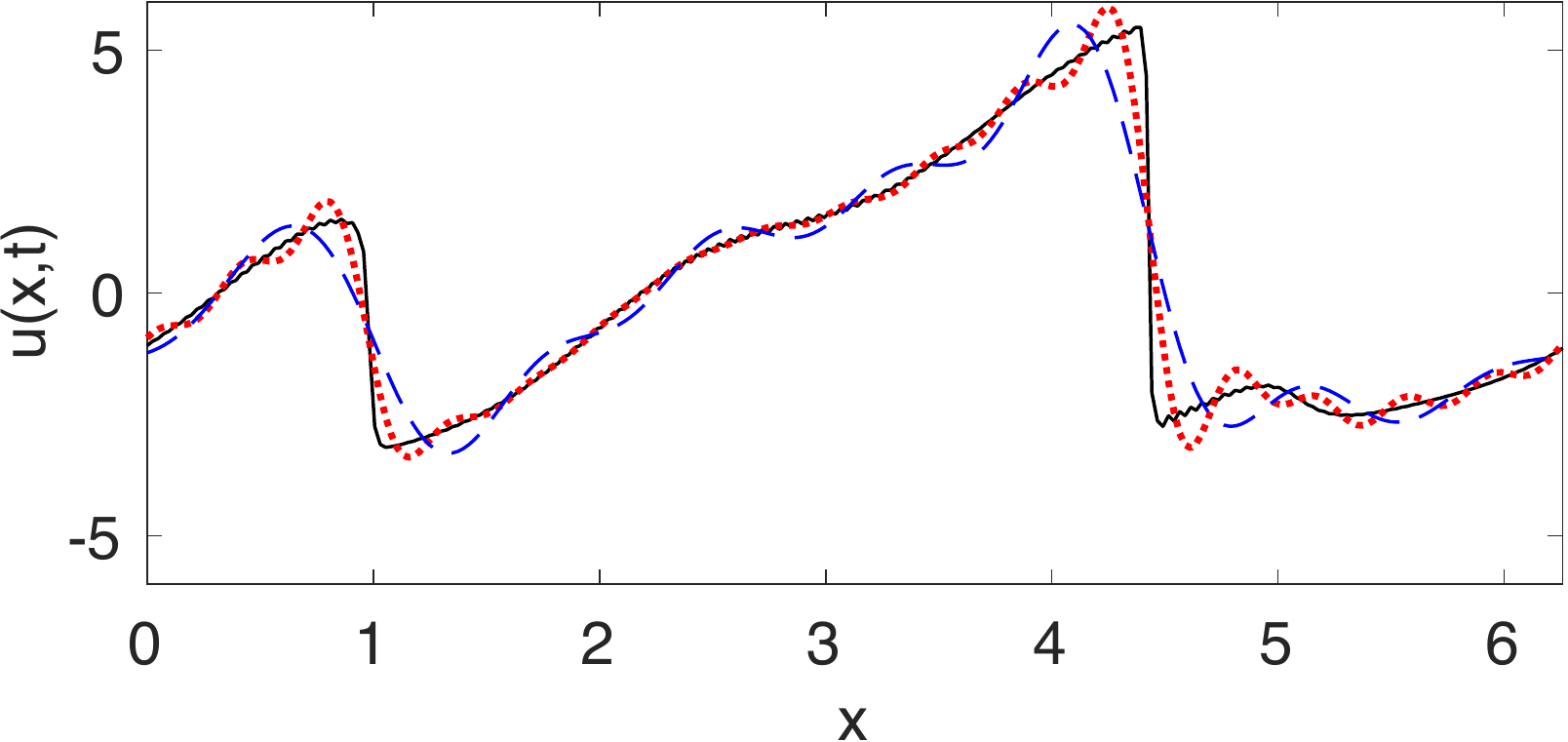} } 
		\subfigure[Stochastic force  $\sigma=0.2$]{\includegraphics[scale=.47]{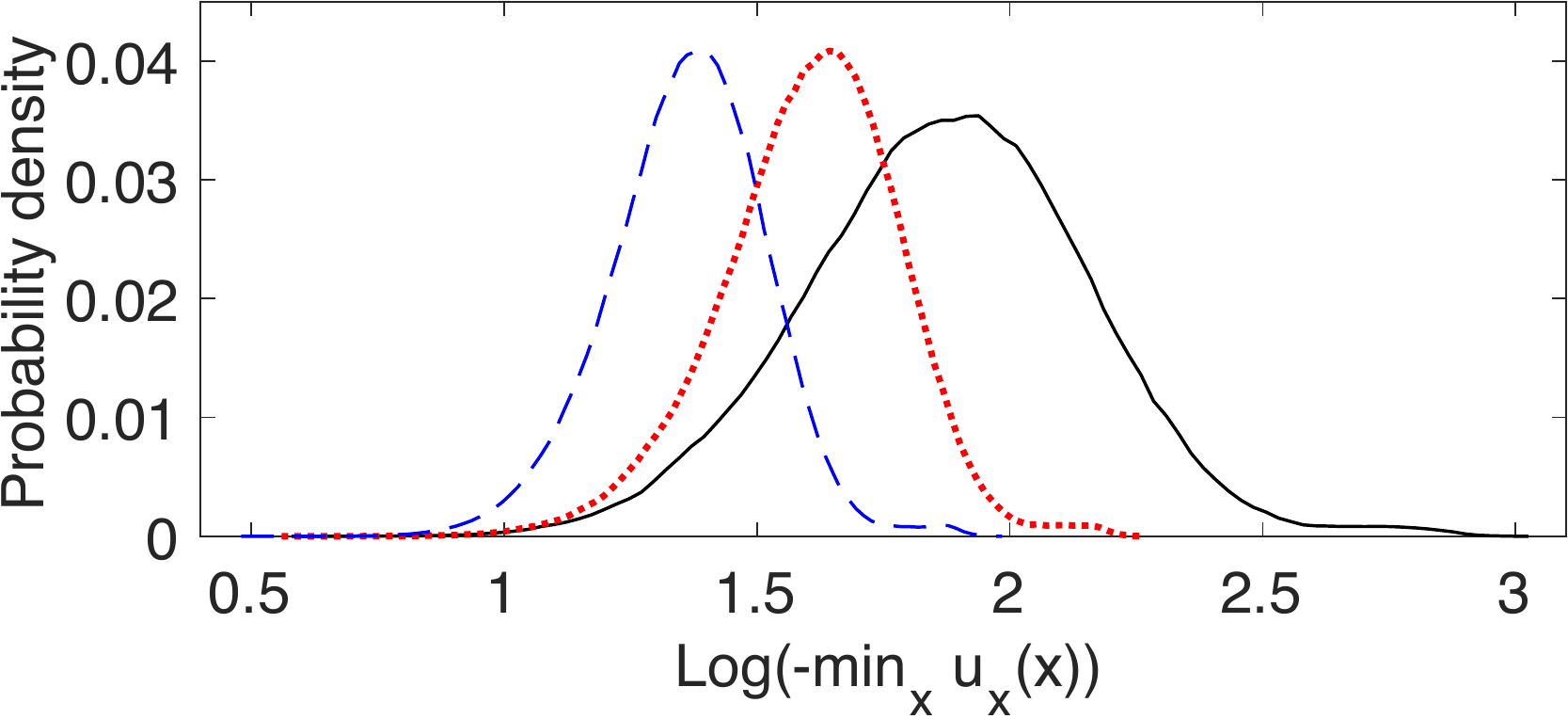}} 
		\subfigure[Stochastic force $\sigma=1$]{\includegraphics[scale=.47]{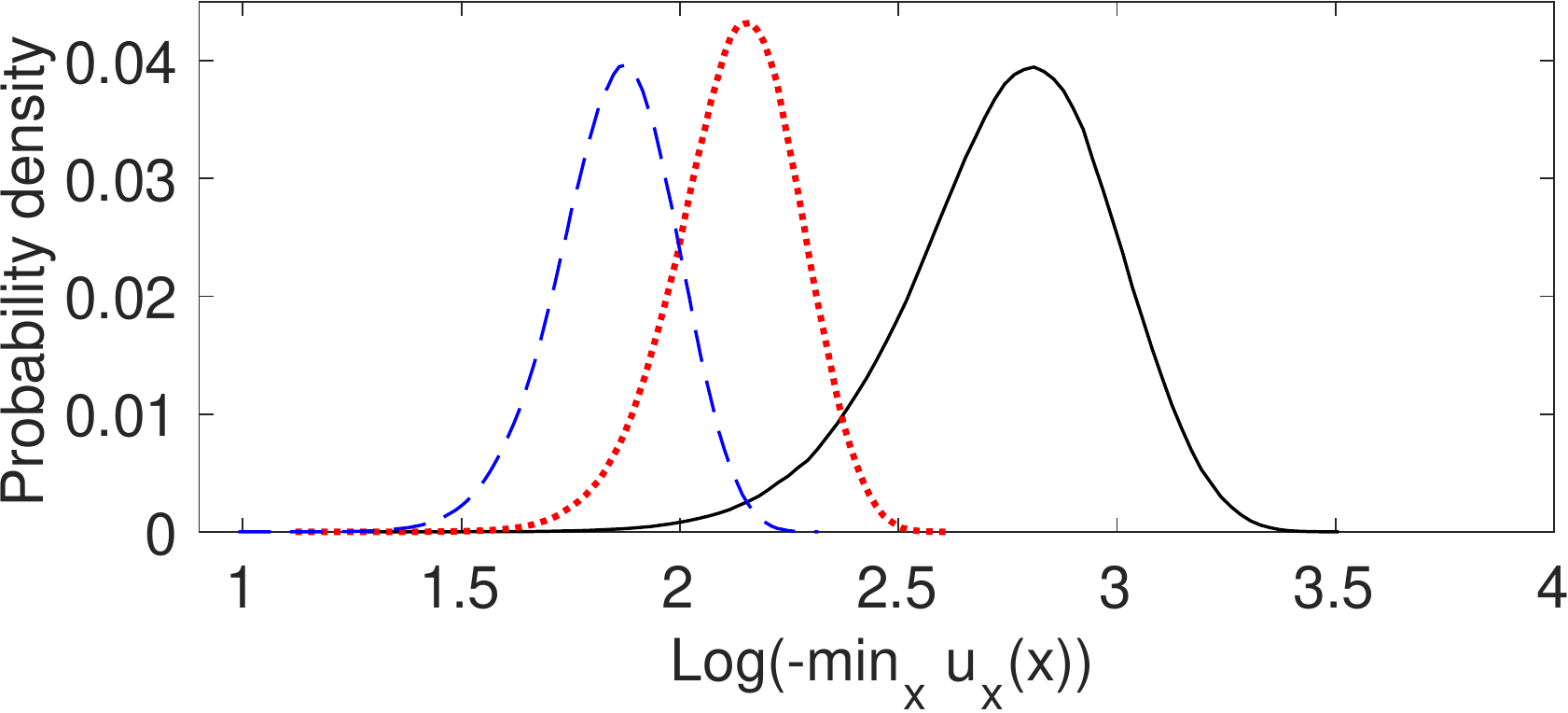} } 
\caption{Typical shocks (Panels a-b) and empirical densities of most negative derivatives (Panels c-d) of the full model solution (denoted by FM), and its $K$-mode and $2K$-mode projections with $K=8$. The full model is simulated with $N=128$ Fourier modes and the viscosity is $\nu=0.02$; see Table~\ref{tab:RMsettings} for the value of other parameters. 
 The density functions have similar shape, but the most negative derivatives of the full model solution are almost a magnitude larger than the other two.  We applied logarithm (with base 10) to the absolute value of the most negative derivatives when computing the probability densities for visualization purpose.}   \label{fig:pdf_Du}
\end{figure}

Figure \ref{fig:pdf_Du} shows the empirical probability distributions of the most negative derivatives $D_k^{(m)}(t_l)$ for the full model solution $u_N$, the K-mode and the 2K-mode projections $u_K$ and  $u_{2K}$ (with $k=N$, $k=K $ and $k=2K$, respectively).  They are computed from $M=200$ trajectories with length $T=100$.  The full model solution's most negative derivatives are about a magnitude larger than the other two. The shocks have derivatives at the scale $-10^3$, on the right end of the distribution of the full model solution, which are out of the reach of the low frequency Fourier modes. Thus, the low frequency Fourier modes alone cannot represent the shocks. But the three density functions have similar shapes, suggesting  connections between the most negative derivatives of the full model solution and those computed from its lower-dimensional projections. Our empirical criterion above exploits these connections to detect the trace of shocks.

The empirical threshold values $\tau_k$ are shown in Table \ref{tab:threshold}. The threshold value increases as the resolution increases in the case of strong stochastic force. But in the case of weak stochastic forcing, due to the large variation of the most negative derivative of the full model (see Figure \ref{fig:pdf_Du}(b)), $\tau_{2K}$ turns out to be more negative than $\tau_N$.
 \begin{table}[H]\centering 
\caption{Threshold values in defining shock traces}\label{tab:threshold}
\begin{tabular}{|c| c  c c|}
 \hline
       Forcing strength  & $\tau_N$  &  $\tau_{2K}$  &  $\tau_{K}$  \\
        \hline
$\sigma =0.2$ &  -21.26   & -26.02 &  -15.77   \\ 
$\sigma =1$ & -320.41 & -98.16&  -52.22\\
  \hline
\end{tabular}

\vspace*{2ex}
{\footnotesize The threshold values provided here are computed for the full model solution ($\tau_N$), and its $2K$-mode ($\tau_{2K}$) and $K$-mode projections ($\tau_K$), where $N=128$ and $K=8$. The system parameters are those listed in Table~\ref{tab:RMsettings}.}
\end{table}

\begin{figure}[tbh!]
    \centering 
  \subfigure[Stochastic force  $\sigma=0.2$]{    \includegraphics[width =0.48\textwidth]{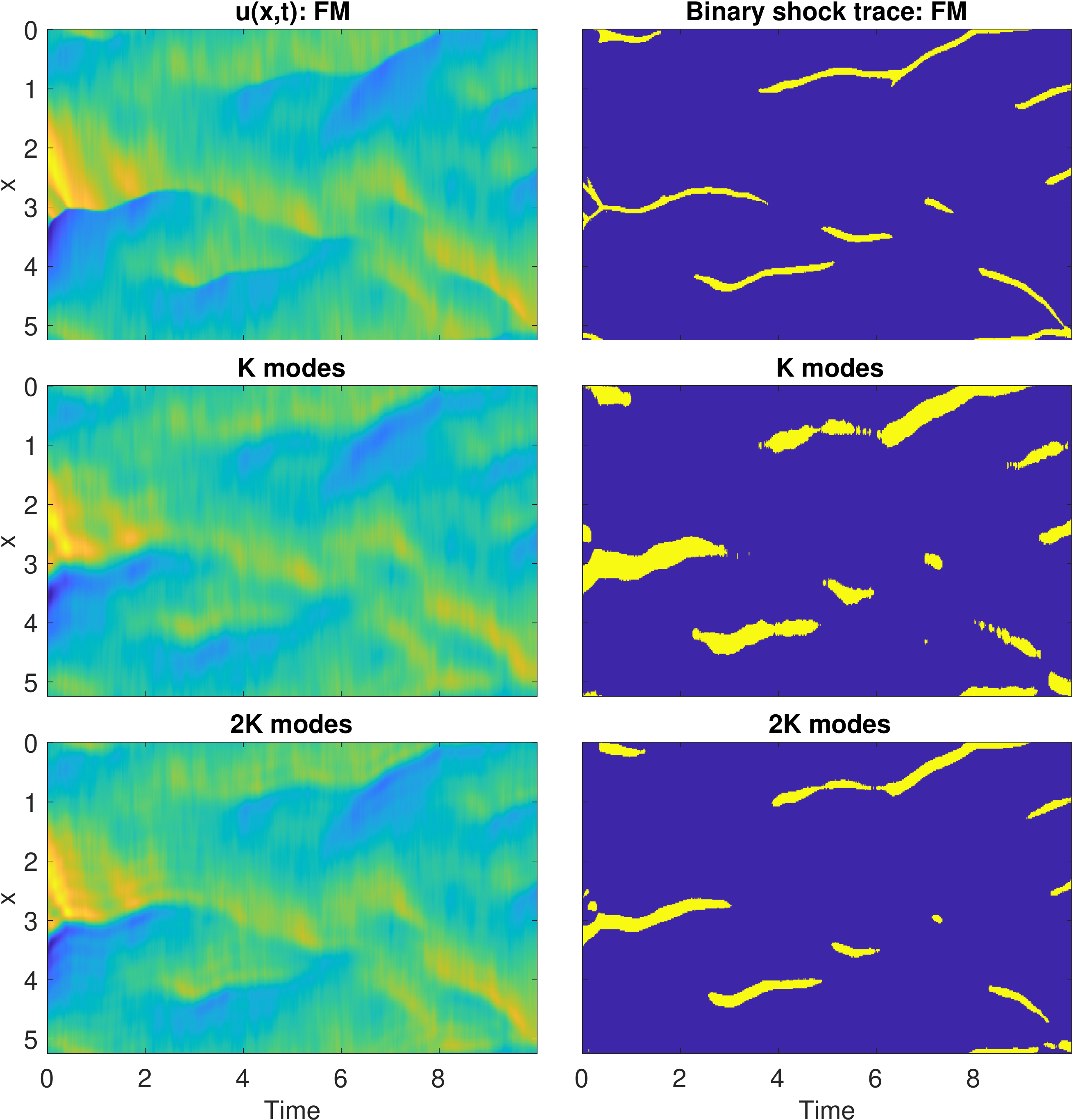} } 
    \subfigure[Stochastic force $\sigma=1$]{      \includegraphics[width =0.48\textwidth]{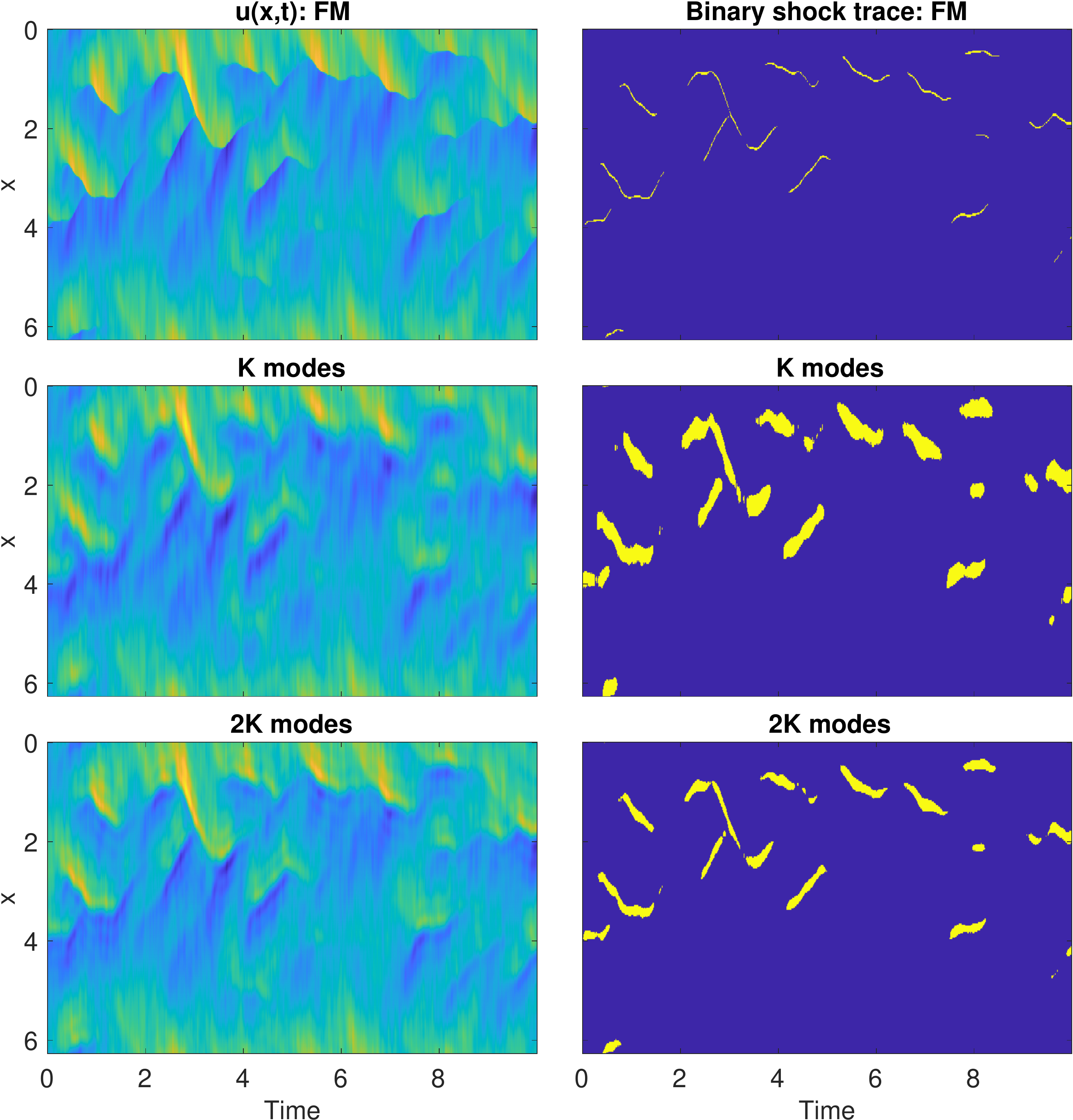} } 
    \vspace{-3mm}
\caption{\small  The full model solution, and its $K$-mode and $2K$-mode projections, in space-time plot  (left column) and their binary shock trace (right column) obtained by their thresholds in \eqref{eq:tau}. The $K$-model projection already encodes the shock trace information, providing a basis for using reduced models to predict the shock trace.
Here $K=8$ and the full model has $N=128$ Fourier modes.} 
\label{fig:traceShock_true}
\end{figure}

With the threshold value $\tau_k$ for a $k$-mode projection $u_k$ as above, we define a \emph{shock trace indicator} function:
\begin{align}\label{eq:binaryShockTrace}
\mathbf{1}_{S_{u_k}}(x,t) =     \left\{
    \begin{array}{ll}
      1, & (x,t)  \in S_{u_k} ,\\[1ex]
      0, &  (x,t) \notin S_{u_k},
    \end{array}
    \right. \, \text{ with }S_{u_k} = \{ (y,s)\in [0,2\pi]\times [0,T]: \partial_x u_k(y,s)<\tau_k\}.
\end{align}
The binary shock trace indicator function has value $1$ in the set $S_{u_k}$, a neighborhood where shock occurs, and is zero otherwise. It is applicable to both the full model solution and its truncation  (see Figure \ref{fig:traceShock_true}, bottom row), as well as the solutions of a reduced model to be studied in later sections.

Figure  \ref{fig:traceShock_true} shows that our thresholds defined in \eqref{eq:tau} can detect the trace of the viscous shocks from either the full model solution or its K-mode and 2K-mode projections. These thresholds, set to be one standard deviation above the corresponding mean, are based on an empirical balance between accuracy of detection and robustness of tolerating the random perturbations from the stochastic force. As a result,  the binary shock trace becomes wider as the number of Fourier modes of the truncated solution decreases from $2K$ to $K$ (with $K=8$ here). A more negative threshold (e.g., one standard deviation below the mean) can lead to narrower shock traces; but for the parameter regimes considered, the corresponding shock traces identified by the K-mode and 2K-mode projections match less good with those identified from the full model solution.

In our numerical tests in Section \ref{sec:numerics}, we will consider reduced models with $K=8$ to show that the proposed closure model is able to faithfully reproduce the shock traces as revealed by the projection of the true solution on the first $K$ modes while a standard $K$-mode Galerkin truncation fails.

\begin{remark} \label{Rmk:generalized_tauk} 
Our empirical threshold provides a preliminary criterion for the detection of shocks for different resolutions. The threshold values are resolution dependent. Here we set the threshold value to be one standard deviation above the mean of the most negative spatial-derivatives as in {\rm\eqref{eq:tau}} for simplicity. One can further refine the definition by multiplying for instance a resolution-dependent factor $\lambda_k$ to the standard deviation, leading to $\tau_k =  \overline{D}_k + \lambda_k \eta_{k}$. One could also set the threshold to be simply the maximal of the corresponding density function such as shown in Panels c-d of Figure~\ref{fig:pdf_Du}.  
\end{remark}

\section{Data-Driven Reduced Model with Closure  Modeling}\label{sec:NAR}
To reduce the computational cost and facilitate an efficient ensemble forecast of the shock trace, a data-driven nonlinear reduced model is developed here that involves the dynamics of only the leading $K$ Fourier modes. The development of such a reduced model is based on a nonlinear autoregression (NAR) modeling framework, which was introduced in \cite{CL15,Lu20Burgers,LinLu21}. The key idea in the construction of the NAR model is a parametric approximation of the discrete-time flow map of the leading $K$ modes. In addition to the reduction of the spatial dimension, a much larger numerical integration time step is used in simulating the NAR model, which enhances the overall computational efficiency.
It is important to note that although the discrete-time flow map is a time-varying infinite-dimensional functional that depends on the $K$ modes and the trajectory of the stochastic forces, our inference of the NAR model does not suffer from the curse of dimensionality because we make use of the nonlinear structure in the equation to derive an informative parametric approximation of the flow map, as will be seen in Section~\ref{sec:NAR_lowModes}. 

\subsection{Discrete-time flow map and its approximation}
We start with the general framework of the model reduction in terms of discrete-time flow map approximation, which can be used to construct the reduced models with suitable closure terms.

To simplify notion, we write the stochastic Burgers equation \eqref{SBE} in an operator form as
\begin{equation} \label{eq_diss} 
\frac{\mathrm{d} u}{\mathrm{d}t} + Au = B(u) + f, \ \  u(0)= u_0,
\end{equation}
where the linear operator $A:H^1_p\to L^2$ and the nonlinear operator $B: H^1_p\to L^2$ are defined by 
\[
A= -\nu \Delta,\quad B(u) = -(u^2)_x/2,
\]
and $f$ is the stochastic force. Here  $H^1_p$ denotes the standard Sobolev space of periodic functions on $[0,2\pi]$ that are square integrable and whose first-order weak derivatives are also square integrable. 

We first decompose the system into the resolved and the unresolved scales (e.g., the low frequency and high frequency Fourier modes). That is, we write the solution in the form
\begin{equation*}
u=Pu+ Qu=v+w,
\end{equation*}
where $P$ and $Q=I-P$ denote the projection operators to the resolved and unresolved scales, respectively. For example, in terms of Fourier expansion $u=\sum_{|k|=1}^\infty \widehat u_k(t) e^{i kx}$, we set $P$ and $Q$ be the projections from $H_p^1$ to the low and high wavenumber Fourier modes, with 
\begin{equation}
v=\sum_{|k|=1}^K \widehat u_k(t) e^{ikx}, \quad  w=\sum_{|k|=K+1}^\infty \widehat u_k(t) e^{i kx}.
\end{equation} 
With these notations, we can write the system \eqref{eq_diss} as
\begin{subequations}\label{eq:full-sys}
 \begin{align}
\frac{\mathrm{d}v}{\mathrm{d}t}&=-PAv +PB(v) + Pf + [PB(v+w)-PB(v)],   \label{lowModes} \\
\frac{\mathrm{d}w}{\mathrm{d}t}&=-QAw+QB(v+w) + Qf.  \label{eq:highModes}
\end{align}
\end{subequations}
The equation \eqref{lowModes} for the resolved variables $v$ is not closed since it depends on $w$ due to nonlinear coupling brought by the nonlinear operator $B$. Thus, a closure model for the resolved-scale variables aims to approximate $[PB(v+w)-PB(v)]$ (see e.g.,\cite{GHM10,Har16,CL15,mou2020dd,chekroun2020variational}) by a function of $v$ and possible additional noise terms. The Mori-Zwanzig formalism for deterministic systems \cite{CH13,CHK02,chorin2003conditional,LinLu21,santos2021reduced} shows that an exact closure model involves memory effects (i.e., the dependence on the history of $v$) and the uncertainties from the unknown initial condition of $w$. Both the memory effects and uncertainties are difficult to model from physical principles. Therefore, date-driven approaches, combined with physical insights, have lead to various constructions of closure models by statistical learning of the discrete-time flow map, either by parametric models \cite{CL15, LLC17, Lu20Burgers} or by nonparametric machine learning representations \cite{ma2018model, harlim2020machine, gouasmiPrioriEstimation2017}. These closure models can describe the statistical and dynamical properties of the resolved scales, leading to accurate prediction of $v$ with uncertainty being quantified.

The NAR reduced model we adopted here aims to approximate the discrete-time flow map of $v$.  More precisely, let $t_n= n\delta$ for $n=1,2, \ldots$ with $\delta$ being one time step. The discrete-time flow map of $v$, according to  \eqref{eq:full-sys}, takes the form
 \begin{align}
v(t_n) & = \mathcal{G}(v(t_{n-1}), w(t_{n-1}), f(s)_{s\in [t_{n-1}, t_n)}), \label{eq:flow_v}  
\end{align}
where $\mathcal{G}$ is a functional forwarding the flow of the resolved variable $v$ from time $t_{n-1}$ to $t_n$, which depends on the current state $(v(t_{n-1}), w(t_{n-1}))$ and the trajectory of the stochastic forcing on the time interval $[t_{n-1}, t_n)$. 
 We seek a parametric function $F$ to approximate the flow map $\mathcal{G}$ that is independent of the unresolved variable $w$, but instead, is dependent on the history of the resolved variable $v$ and the stochastic force. For this purpose, we denote by $v^{1:n-1}$ the discrete-time trajectory of $v$ at $t_1,\ldots, t_{n-1}$, and by $f^{1:n}$ the discrete-time trajectory of $f$ at $t_1,\ldots, t_{n}$. Denote also $v^n = v(t_n)$. We aim to construct $F$ such that 
 \begin{align}
v^n& \approx F(\theta,v^{1:n-1}, f^{1:n}) + g^n , \label{eq:dis_v} 
\end{align}
where $\theta$ is a multivariate parameter to be estimated from data and the precise form of $F$ is inspired by a Picard approximation of the unresolved variable $w$ as detailed in Section~\ref{sec:NAR_lowModes} below. The additional process $\{g^n\}$ aims to represent the residual $\mathcal{G}(v(t_{n-1}), w(t_{n-1}), f(s)_{s\in [t_{n-1}, t_n)}) - F(\theta,v^{1:n-1}, f^{1:n})$. In this paper, $g^n$ is set to be Gaussian for simplicity, but other forms are also possible, taking the moving average \cite{CL15,LLC16} or stochastic Stuart-Landau oscillators \cite{chekroun2021stochastic} as examples.

 \subsection{The parametric reduced model with closure modeling terms}\label{sec:NAR_lowModes}

The construction of a closure model for $v$ of the form \eqref{eq:dis_v} is through a statistical learning from data, which was originally introduced in \cite{Lu20Burgers}. The construction consists of three steps:
\begin{enumerate}
  \item deriving a family of parametric functions from the Picard approximation of the high frequency modes,
  \item estimating the parameters by maximizing the likelihood of the data of $v$, and
  \item selecting the model that best fits the data.
\end{enumerate}

\subsubsection{Derivation of the parametric model}
We first derive the parametric function $F(\theta,v^{1:n-1}, f^{1:n})$ in \eqref{eq:dis_v} from numerical integrators of the full model. In view of Eq.\eqref{lowModes}, a closure model has to represent the residual $PB(v+w)-PB(v)$ in term of $v$ and its history.  Note that $w(t)$ is a functional of the history of $v$ by integrating \eqref{eq:highModes}: 
\begin{equation}\label{eq:w_intg}
\begin{array}{ll}
w(t) &= e^{-QAr} w(t-r) + \int_{t-r}^t e^{-QA(t-s)} [QB(v(s)+w(s))+ Qf(s)]ds, 
\end{array}
\end{equation}
where $r\in [0,t]$. Given $w(t-r)$ and a trajectory $(v(s), Qf(s), s\in [t-r, t])$,  the Picard iteration provides us an explicit approximation of  $w(t)$ as a functional of the trajectory of $v$. That is, the sequence of functions $\{w^{(l)}\}$, defined by
\begin{equation} \label{eq:w_picard}
w^{(l+1)}(t) = e^{-QAr} w^{(l)}(t-r) + \int_{t-r}^t e^{-QA(t-s)} [QB(v(s)+w^{(l)}(s))+ Qf(s)]ds,
\end{equation}
with  $w^{(0)}(s)=0$ for $s\in [t-r, t]$, converges under suitable conditions to the function $(w(s), s\in [t-r, t])$ as $l\to \infty$. In particular, the first Picard iteration provides a closed representation 
\begin{equation} \label{Eq_LP_integral}
w^{(1)}(t) = \int_{t-r}^t e^{-QA(t-s)} [QB(v(s))+ Qf(s)]ds.
\end{equation}
While a rigorous convergence analysis of the above Picard iteration is beyond the scope of the current paper, we will show below that the NAR model built from the first iteration \eqref{Eq_LP_integral} can already provide good performances. Note that integrals of the form \eqref{Eq_LP_integral} also arise naturally in the approximation of center manifolds as mentioned in Remark~\ref{Rmk:CM} below. 

Now, substituting $w$ in \eqref{lowModes} by $w^{(1)}$, we obtain an approximate closed  integro-differential equation for $v$. 
We then discretize this integro-differential equation in a parametric way to obtain a parametric function that aims to approximate the discrete-time flow map of $v$. For this purpose, let $\delta = t_n-t_{n-1}$ denote the time-step size and let $r = p\delta$. 
By parametrizing a Riemann sum approximation of the above integral \eqref{Eq_LP_integral} for $w^{(1)}(t)$, we get 
$w^{(1)}(t_{n})  \approx  \sum_{j=0}^p c_j e^{-QA j \delta} [ QB(v(t_{n-j})) +Qf(t_{n-j})]$. We denote
\begin{equation}\label{eq:w_intgAppr}
\begin{aligned}
w^{(1),n} &=  \sum_{j=0}^p c_j e^{-QA j \delta} [ QB(v^{n-j}) +Qf^{n-j}],
\end{aligned}
\end{equation}
where $c_j\in \R$ are parameters to be estimated from data.
Furthermore, to keep linear dependence on the parameters,  we approximate the term $PB(v+w)-PB(v)$ in \eqref{lowModes} as follows by ignoring the nonlinear self-interaction term $P w\partial_{x}w$:
\[PB(v+w^{(1)})-PB(v) \approx - P( \partial_{x} v w^{(1)}  + v\partial_{x} w^{(1)}).
\]
Then, by parametrizing a numerical integrator of \eqref{lowModes}, we obtain:
\begin{subequations}
 \begin{align} \label{eq:dis-v_para}
v^n&=v^{n-1}+ \delta[R^{\delta}(v^{n-1} ) + Pf^{n-1} +\Phi^{n-1} ] + g^n,   \\
\Phi^{n-1}&=  \partial_{x} v^{n-1}  w^{(1),n-1}+v^{n-1} \partial_{x} w^{(1),n-1}   + \sum_{j=1}^{p} [c^v_j v^{n-j}+c_j^R R^{\delta}(v^{n-j} ) + c^f_j Pf^{n-j}],\label{eq:QBvLag}
\end{align}
\end{subequations}
where the nonlinear function $R^{\delta}(\cdot)$ comes from a numerical integration of the deterministic truncated Galerkin equation $\frac{\mathrm{d}v}{\mathrm{d}t} \approx -PAv +PB(v)$ at time $t_{n-1}$ and with time-step size $\delta$. Here the extra parametrized terms $[c_j^R R^{\delta}(v^{n-j} ) + c^f_j Pf^{n-j}]$ aim to further account for the memory. Note that the high frequency component $Qf$ of the noise does not enter the closure of the low modes because $PQf=0$ in the linear approximation.

Next, we rewrite \eqref{eq:dis-v_para} and \eqref{eq:QBvLag} in terms of Fourier modes, and we further parametrize the modes component-wisely to account for the different dynamics between modes.
 Denote $\wv^n = (\wv^n_k, |k|\leq K)\in \C^{2K}$ the low-modes in the reduced model that approximates the original low modes $(\hatu_k (t_n), |k|\leq K)$. The reduced model is in the form of a nonlinear autoregression (NAR) model, which reads
 \begin{subequations} \label{eq:NAR}
 \begin{align}
\wv^n_k & = \wv^{n-1}_k+  \delta [ \widehat R^{\delta}_k(\wv^{n-1}) +   \widehat{f}^{n-1}_k + \Phi^{n-1}_k ]+ \widehat{g}^n_k, \quad 1\leq k\leq K, \label{eq:NAR_v} \\
\Phi^{n-1}_k & =  \sum_{j=1}^{p}\left[ c^v_{k,j}\wv^{n-j}_k + c^R_{k,j} \widehat R^{\delta}_k( \wv^{n-j}) + c^f_{k,j} \widehat{f}^{n-j}_k+  c^w_{k,j} \sum_{\substack{ |k-l|\leq K, K< |l| \leq 2K \\ \text{ or }  |l|\leq K, K< |k-l| \leq 2K} } {\widetilde v^{n-1}_l \widetilde v^{n-j}_{k-l}} \right], \label{eq:NAR_linearQ}
\end{align}
\end{subequations}
 where $  \widehat R^{\delta} $ comes from the ETD-RK4 integrator of the $K$-mode truncated system, $\widehat{f}^n_k$ and $\widehat{g}^n_k$ denote the $k$-th mode coefficient in the Fourier transform of respectively $f^n$ and $g^n$, and the notation $\widetilde v^{n-j}_l$  is defined by  \begin{equation}  \label{eq:ks-ansatz-d}
  \widetilde{v}^{n-j}_{k}=
    \left\{
    \begin{array}{ll}
      \wv^{n-j}_k~, & 1\leq k\leq K;\\[1ex]
      \frac{i k}{2}e^{-\nu k^2 j\delta} \sum_{|l|\leq K, |k-l|\leq K}\wv^{n-j}_{k-l} \wv^{n-j}_{l} , & K < k \leq 2K.
    \end{array}
    \right.
  \end{equation}
We set $\wv^n_{-k}= (\wv^n_k)^*$ (with the sup-script $^*$ denoting
complex conjugate).  Here the coefficients vary component-wisely to allow more flexibility in fitting data.

The equation \eqref{eq:NAR} defines a nonlinear autoregression type model when the residual term $\{g^n\}$ is modeled by independent Gaussian noise.  One may further improve the model in two directions: (1) include more nonlinear terms by parametrizing higher order Picard iteration than the first iteration; (2) consider spatial correlation between the components of $g$
or by using moving average models \cite{CL15,LLC16}.   We assume for simplicity that the $g$ has independent components, so that the coefficients can be estimated by the least squares fitting.

\begin{remark} \label{Rmk:CM}
The parameterization of the unresolved variable $w$ given by \eqref{Eq_LP_integral} is closely related to the Lyapunov-Perron integrals arising from the approximation of center manifolds {\rm\cite{Hen81,MW05,CLW15_vol1}}, in which the parameter $r$ is pushed to infinite and the true dynamics of $v$ involved in the integral \eqref{Eq_LP_integral} is replaced by solutions from a linearized equation; see e.g.~{\rm \cite[Theorem 1]{chekroun2020variational}} and {\rm \cite[Theorem 6.1]{CLW15_vol1}}. Treating $r$ in \eqref{Eq_LP_integral} as a free parameter to be optimized using solution data leads to deformations of the manifold that aims to approximate the unresolved dynamics in an optimal way. Such a data-driven optimization is the key of the parameterizing manifold approach proposed in {\rm \cite{CLW15_vol1,chekroun2020variational}} for dimension reduction of stochastic and chaotic systems. The difference here is that we only utilize the functional form of the parameterization, but optimize all of the associated coefficients using data.  
\end{remark}

\subsubsection{Parameter estimation} \label{Sec_param_est}
The coefficients in \eqref{eq:NAR} are estimated by maximizing the likelihood of data.  The data can be either a long trajectory or many independent short trajectories. We denote the data consisting of $M$ independent trajectories by
\begin{equation}\label{eq:data}
\text{Data: } \quad
 \{\wv^{1:N_t,m}, \widehat f^{1:N_t,m}\}_{m=1}^{M}  = \{ (\hatu_k^{(m)}(t_{1:N_t})), (\widehat{f}_k^{(m)}(t_{1:N_t})), k=1,\ldots, K \}_{m=1}^{M}, 
\end{equation}
 where $m$ indexes the trajectories and $N_t$ denotes the number of steps for each trajectory.

The maximal likelihood estimator (MLE) of the coefficients is computed by least squares since the model depends linearly on them \cite{Lu20Burgers}. More precisely, writing \eqref{eq:NAR_linearQ} as 
\[
\Phi_k^n(\boldsymbol\theta_k) = \sum_{j=1}^{4p}\boldsymbol\theta_{k,j} \boldsymbol \Phi^n_{k,j}
\]
 with $\boldsymbol\theta_k =(c^v_{k,j}, c^R_{k,j}, c^f_{k,j},c^w_{k,j}, j=1,\ldots, p)\in \R^{4p}$ and $\boldsymbol\Phi^n_{k} = (\wv^{n-j}_k, \widehat R^{\delta}_k( \wv^{n-j}),\widehat  f^{n-j}_k, \sum_{k,l} {\widetilde v^{n-1}_l \widetilde v^{n-j}_{k-l}} , j=1,\ldots, p) \in \C^{4p}$ denoting the parametric terms, we compute the MLE as
\begin{equation}\label{eq:MLE}
\begin{aligned}
\widehat {\boldsymbol\theta}_k &= (\bA_k)^{-1}\bb_k, \quad 1\leq k \leq K, \\
\widehat \sigma^g_{k} & = \frac{1}{MT} \sum_{n,m=1}^{T,M} \|\wv_k^{n,m} - \big(\wv^{n-1,m}_k+   \delta  \widehat R^{\delta}_k(\wv^{n-1,m}_k) + \delta\widehat f^{n,m}_k +  \delta \boldsymbol\Phi^{n,m}_k (\widehat{\boldsymbol \theta}) \big)\|^2
\end{aligned}
\end{equation}
where the normal matrix $\bA_k$ and vector $\bb_k$ are defined by
\begin{align}
\bA_k(j',j) & =   \frac{\delta}{MT} \sum_{n,m=1}^{T,M} 	\langle \boldsymbol\Phi^{n,m}_{k,j'}, \boldsymbol\Phi^{n,m}_{k,j}\rangle, \quad 1\leq j', j \leq 4p,  \label{eq:normalMat}\\
\bb_k(j) & =    \frac{1}{MT} \sum_{n,m=1}^{T,M}	\langle \wv^{n,m}_k - \wv^{n-1,m}_k+   \delta \widehat R^{\delta}_k(\wv^{n-1,m}_k) +\delta\widehat f^{n,m}_k , \boldsymbol\Phi^{n,m}_{k,j}\rangle.  \notag
\end{align}
In practice, we use pseudo-inverse or regularization to solve the least squares problem in \eqref{eq:MLE} when $\mathbf A_k$ is ill-conditioned. 

By fitting parameters to the data of true solution, we obtain the optimal function in the parametric family, correcting the numerical error and model error.

\subsubsection{Model selection} \label{Sec_model_selection}
The model selection step aims to determine the time lag $p$ and remove the redundant terms in the model in \eqref{eq:NAR_linearQ}.  We select the simplest model that  fits the data the best in the sense that: (i) it reproduces the statistics such as the energy spectrum, the marginal invariant densities and temporal correlations; (ii) the estimator converges as the data size increases.

For the settings in Table \ref{tab:RMsettings}, we take $p=1$ to yield the simplest models, i.e., our reduced models are Markovian. Non-Markovian models can improve the results and we refer to \cite{Lu20Burgers} for further discussions on the memory length and parameter convergence. With $p=1$, the model selection step suggests that we estimate the parameters  $(c^{v}_{k,1}, \, c^{R}_{k,1}, \, c^f_{k,1},\, c^{w}_{k,1})$ for $1\leq k \leq K$, along with the variance of the residuals of the $K$ modes. The estimators converge fast as the data size increases, either in the number of trajectories or in the length of a trajectory.

The estimators from data consisting of 512 trajectories, each with length 160 time units and with an initial condition down-sampled from a long trajectory, are shown in Table \ref{tab:parNAR_sigma_p2} and \ref{tab:parNAR_sigma1}. All the estimators depend on the strength of the stochastic force, but they have the following common features. The estimators  $(\widehat {c}^v_{k,1}, k=1,\ldots, K)$ are all negative and are more negative as $k$ increases, so that the corresponding linear terms all dissipate energy. The estimator $(\widehat {c}^R_{k,1}, k=1,\ldots, K)$ are round 1. These two estimators, together with the estimators $(\widehat {c}^f_{k,1}, k=1,\ldots, K)$, which are small, act as a calibration of the Euler-Maruyama type approximation in \eqref{eq:NAR}. At last, the estimator $(\widehat {c}^w_{k,1}, k=1,\ldots, K)$ are at a smaller magnitude, but they are important for the NAR model to account for the effects of the unresolved scales. These estimated parameters specify our NAR model to predict the shock trace in the next section. 

\section{Shock Trace Prediction by the Reduced Model} \label{sec:numerics}

In this section, we demonstrate the skill of the constructed NAR reduced model \eqref{eq:NAR} for predicting the binary shock trace in the spatiotemporal solution fields of the stochastic Burgers equation \eqref{SBE}. The accuracy is quantified by computing the rate of false positive and false negative events as detailed in Section~\ref{Sec_FP_FN}. We consider two prediction scenarios, one with noise-free initial data (Section~\ref{Sec_noise_free_IC}) and the other with noisy data, for which we use data assimilation techniques to estimate the state and make ensemble prediction to quantify the uncertainty  (Section~\ref{Sec_noisy_IC}).  

Our choice of the reference shock traces is based on the assessments carried out in Section~\ref{sec:shockTrace}. It was shown there that the binary shock traces computed from the $8$-mode projection of the true solution already provide decent indication of the time and locations for the occurrence of shocks in the true solution for both of the two forcing scenarios $\sigma=0.2$ and $\sigma = 1$ within the parameter setup given by Table~\ref{tab:RMsettings}; see the second row in Figure~\ref{fig:traceShock_true}. For this reason, we take the binary shock traces computed from the $8$-mode projection of the true solution as the ground truth; and we set the first $8$ Fourier modes as the resolved modes for the NAR reduced model. A truncated Galerkin system of \eqref{DFM} for the same resolved modes, referred to as the truncated system below, is used as a comparison. As will be shown in Section~\ref{Sec_noise_free_IC} and Section~\ref{Sec_noisy_IC}, the NAR reduced model consistently outperforms the truncated system for both the weak noise forcing case ($\sigma = 0.2$) and the strong noise forcing case ($\sigma = 1$).  

Besides visual comparison of the binary shock trace plots, we provide below a more quantitative way of assessing the performance, which offers statistics that reveals the robustness of the good performance achieved by the proposed NAR reduction framework.

\subsection{Quantification of shock trace prediction skills by false positive/negative rate} \label{Sec_FP_FN}
We define in this section the false positive and false negative rates in predicting the binary shock traces for a reduced model. For this purpose, we take the binary shock traces computed from a $K$-mode projection of the true solution as the ground truth (with $K=8$ here). We get a \emph{false negative} from a reduced model at a space-time location $(x,t)$ when its predicted value is 0 while the truth is 1, and similarly a \emph{false positive} occurs when its predicted value is 1 while the truth is 0. The \emph{false negative rate} (resp.~\emph{false positive rate}) is the ratio between the number of false negatives (resp.~false positives) and the number of true ones.

More precisely, denote by $u_K(x,t)= \sum_{|k|\leq K} \hatu_k(t) e^{ikx}$ the $K$-mode projection of the true solution, and by $v$ the prediction from a reduced model subject to the same realization of the stochastic force used in the full model. Let $\tau_K$ be the binary shock trace threshold associated with $u_K$ defined by \eqref{eq:tau}. Following the definition of the binary shock trace in \eqref{eq:binaryShockTrace}, we denote the sets $S_{u_K} =\{(x,t)\in  [0,2\pi ]\times [0,T]: \partial_x u_K(x,t) < \tau_K\}$ and $S_v =\{(x,t)\in  [0,2\pi ]\times [0,T]: \partial_x v(x,t) < \tau_K\}$. Then, the false negative rate and false positive rate for the prediction $v$ are defined to be 
  \begin{align}\label{eq:rates}
R_{\text{FN}} = \frac{|S_{u_K} \setminus (S_v\cap S_{u_K})|}{|S_{u_K}|}, \quad   R_{\text{FP}} = \frac{|S_v \setminus (S_v \cap S_{u_K})|}{|S_{u_K}|},
  \end{align}
  where $|\cdot|$ denotes the area of a given set. When these rates are computed numerically, the involved sets and the related areas are approximated through the discrete space-time mesh used for the simulation of the reduced model. 
  
The false negative (FN) rate $R_{\text{FN}}$ simply quantifies the percentage of shock traces in $u_K$ that is missed by the reduced model; and the false positive (FP) rate $R_{\text{FP}}$ quantifies the relative number of fake shock traces predicted by the reduced model with respect to the total number of shock traces in $u_K$. Apparently, an accurate prediction would lead to low rates for both false positives and false negatives. Note that these rates depend on the realization of the stochastic force $f$ used as well as the additional stochastic closure terms (if any) in the reduced model. To account for the uncertainty from the randomness, we will compute these rates from many simulations and present them using box charts in Sections~\ref{Sec_noise_free_IC} and \ref{Sec_noisy_IC} below.
 

\subsection{Prediction from noiseless observations} \label{Sec_noise_free_IC}

In this section, we examine the performance of the NAR reduced model in predicting the shock traces when noiseless observations are used as initial data for the prediction. The case with noisy observation data will be dealt with in Section~\ref{Sec_noisy_IC} below. 

The numerical setup is as follows. The parameter regimes are those given by Table~\ref{tab:RMsettings}. Both the NAR model and the Galerkin truncated system aim to model the dynamics of the first $8$ Fourier modes, and the time-step $\delta$ is set to be 10 times bigger than that for the full model. The projection of the initial data for the full model onto the first $8$ modes are used as the initial conditions for the reduced models, and the stochastic force is the same as that used for the full model after coarsening to the larger time-step $\delta$. Thus, each reduced model produces an approximation to the flow map \eqref{eq:flow_v} for the first $8$ Fourier modes, with time-step $\delta$. 

The NAR model is of the form \eqref{eq:NAR}, with the involved parameters trained according to the procedure described in Section~\ref{Sec_param_est}. We have also set $p=1$ in \eqref{eq:NAR}. Thus no memory terms are involved here (see Section~\ref{Sec_model_selection} for further discussions). The additional noise term $g^n$ in \eqref{eq:NAR} is modeled simply by a Gaussian fitting to the residual. In terms of computational cost, aside from the offline training stage for calibrating the involved parameters, the online simulation time used by the NAR model is only slightly longer than that of the truncated model due to the additional closure terms that are absent in the truncated system. But the NAR model's computational cost  remains orders of magnitude smaller than the full system due to its lower space dimension and larger time-step size. 

\begin{figure}[htb]
    \centering
		\subfigure[Stochastic force  $\sigma=0.2$]{\includegraphics[width =0.48\textwidth]{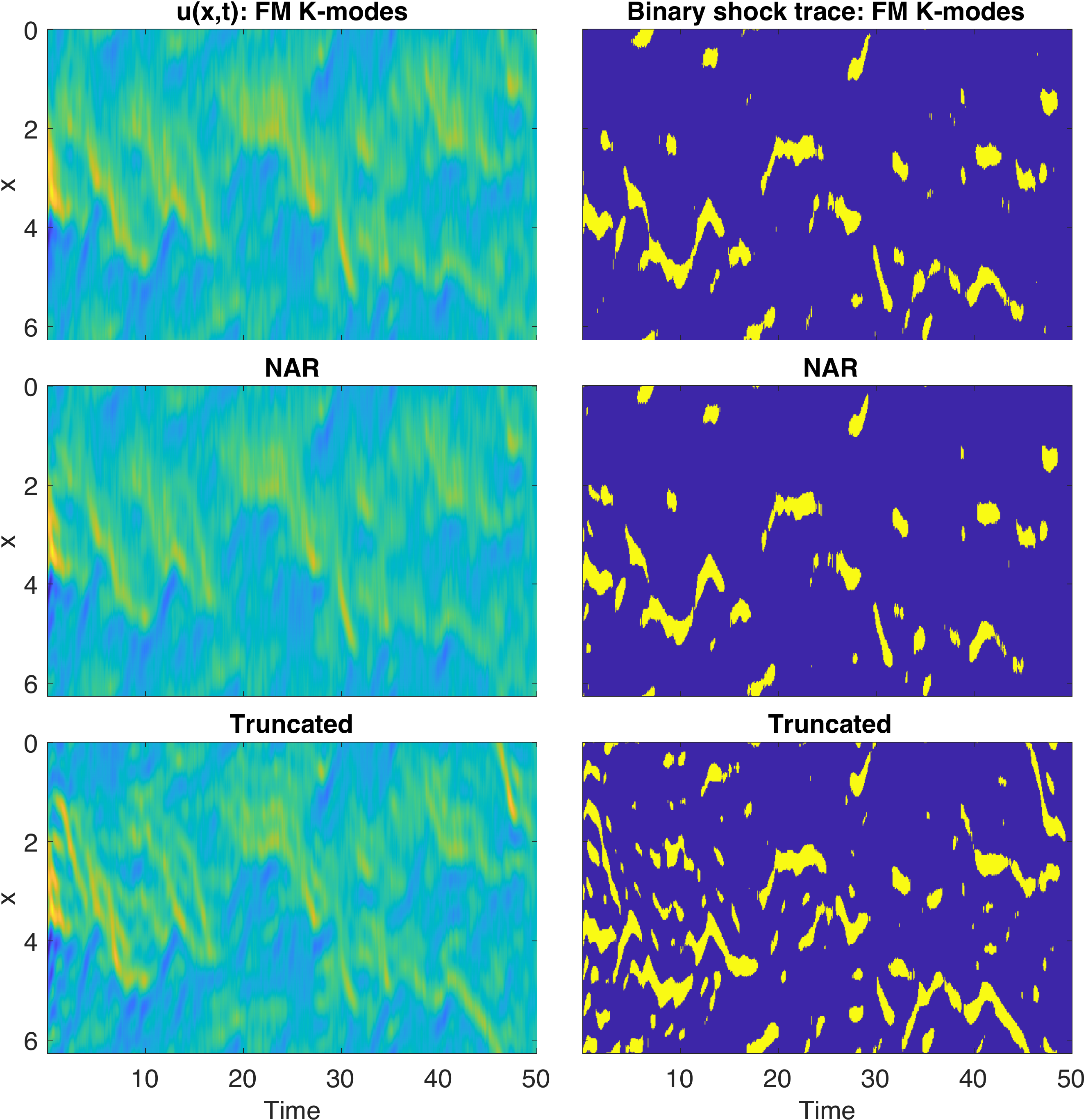}} \hspace{-2mm}
		\subfigure[Stochastic force $\sigma=1$]{\includegraphics[width =0.48\textwidth]{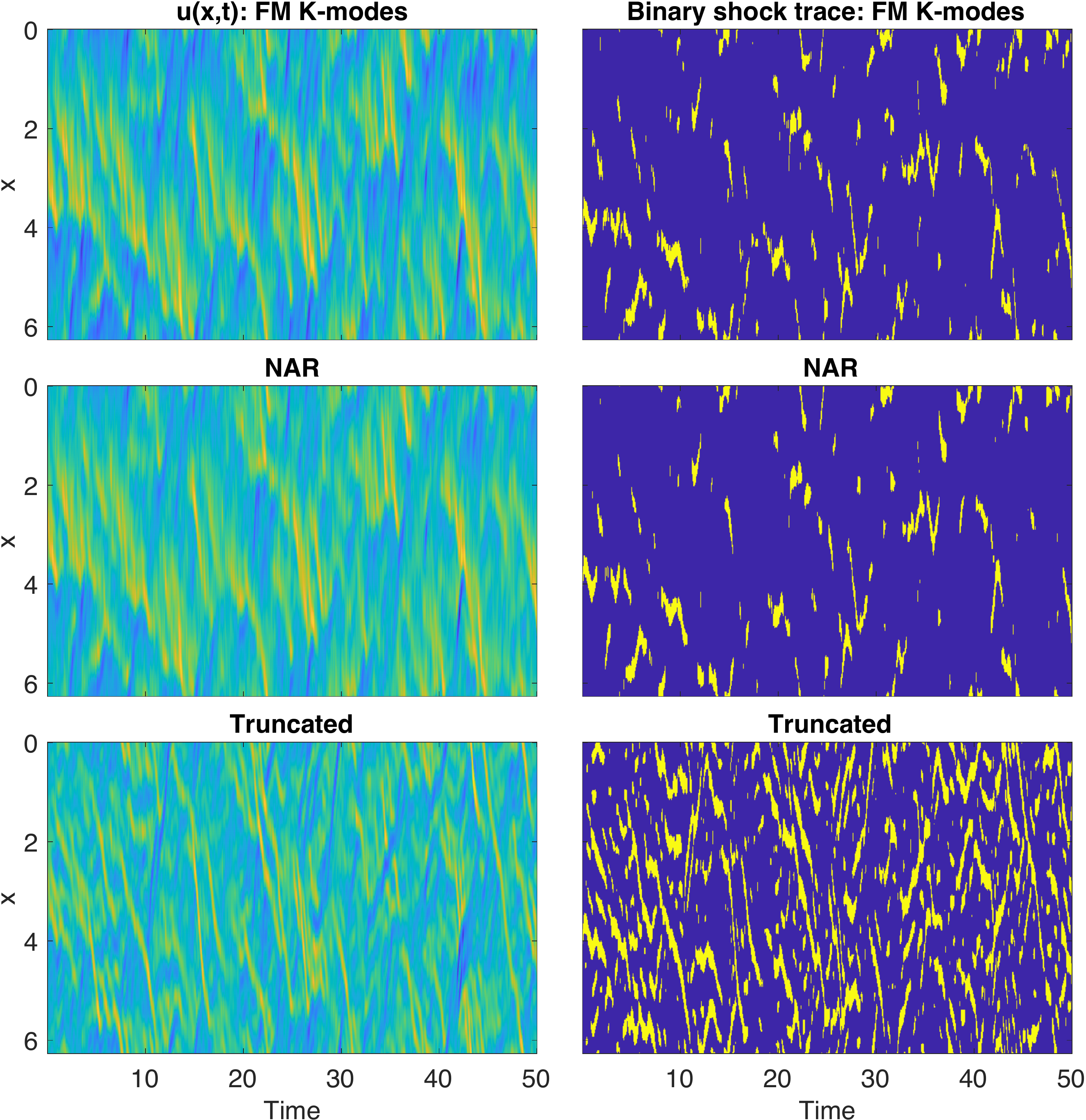}}
		\vspace{-3mm}
\caption{\small Spatiotemporal fields and the corresponding binary traces of viscous shocks from the NAR reduced model (middle row) and the truncated system (bottom row), in comparison with those from the true $K$-mode projection of the solution with $K=8$ (top row). The results are presented for an arbitrarily fixed realization of the stochastic force.  The rates of false predictions of the binary shock traces are shown in Table \ref{tab:rates1}.
 }\label{fig:shock_pred}
\end{figure}

In Figure~\ref{fig:shock_pred}, we present the spatiotemporal fields obtained from the two reduced models and compare them with the  8-mode projection of the true solution field. Also shown are their binary shock traces computed according to \eqref{eq:binaryShockTrace} with $u_{k}$ therein taken to be the spatiotemporal field from either the 8-mode projection of $u$ (top row), or the NAR reduced model (middle row), or the truncated model (bottom row). The threshold $\tau_k$ in \eqref{eq:binaryShockTrace} is taken here to be the $\tau_K$ given by Table~\ref{tab:threshold}. We also report in Table~\ref{tab:rates1} at a more quantitative level the visual results given in Figure~\ref{fig:shock_pred}, using the false negative/positive rates defined in Section~\ref{Sec_FP_FN}. 

\begin{table}[bth!]
\centering
\caption{Rates of false predictions}\label{tab:rates1}
\begin{tabular}{|c| cc | cc|}
 \hline
        &  \multicolumn{2}{c|}{Weak forcing $\sigma =0.2$} &\multicolumn{2}{c|}{Strong forcing $\sigma =1$}  \\
Rate  &  False Positive & False Negative  &  False Positive & False Negative \\ \hline
NAR model & 0.08 & 0.17 &  0.11 & 0.24   \\     
 Truncated system  & 1.06 & 0.31 & 2.71 & 0.55  \\ 
  \hline
\end{tabular}
\vspace{1ex} \\
{\footnotesize The reported rates are computed according to \eqref{eq:rates} for the parameter regimes listed in Table~\ref{tab:RMsettings} with an arbitrarily fixed noise realization. See Figure \ref{fig:shock_pred} for the corresponding binary shock trace plots.}  
\end{table}

For both the weak force regime ($\sigma=0.2$) and the strong force regime ($\sigma=1$), the spatiotemporal field and the binary shock traces of the projected true dynamics are very well reproduced by the NAR model. In contrast, the truncated system already performs visibly less good for the weak force regime, and its predictive skill dramatically decreases in the strong forcing regime. The poor performance of the truncated system is because it does not include any closure term to account for the impact from the unresolved modes. Such closure terms become increasingly important as the forcing strengthens. Indeed, for both forcing scenarios considered, the unresolved modes still retain a significant amount of energy, especially for the strong force regime; see Table~\ref{tab:energy_frac}. 

\begin{table}[bth!]
\centering
\caption{Fraction of energy in the unresolved modes}\label{tab:energy_frac}
\begin{tabular}{|c|c|c|c|}
 \hline
Forcing strength & Mean & Standard deviation & Outliers  \\ \hline
$\sigma = 0.2$  & 3.26\% & 2.04\% & above $10\%$  \\
$\sigma = 1$   &   6.09\% & 3.00\% & above $20\%$  \\      
  \hline
\end{tabular}
\vspace{1ex} \\
{\footnotesize These numbers are computed based on a typical realization of the full model solution over a 400 time-unit window. 
The fraction of energy in the unresolved modes at each time $t$ is computed via $\|u(\cdot, t) - u_K(\cdot, t)\|^2/\|u(\cdot, t)\|^2\times 100\%$, where $\|\cdot\|$ denotes the $L^2$-norm and $K=8$.}    
\end{table}

\begin{figure}[htb]
    \centering
		\subfigure[Stochastic force  $\sigma=0.2$]{\includegraphics[width =0.45\textwidth]{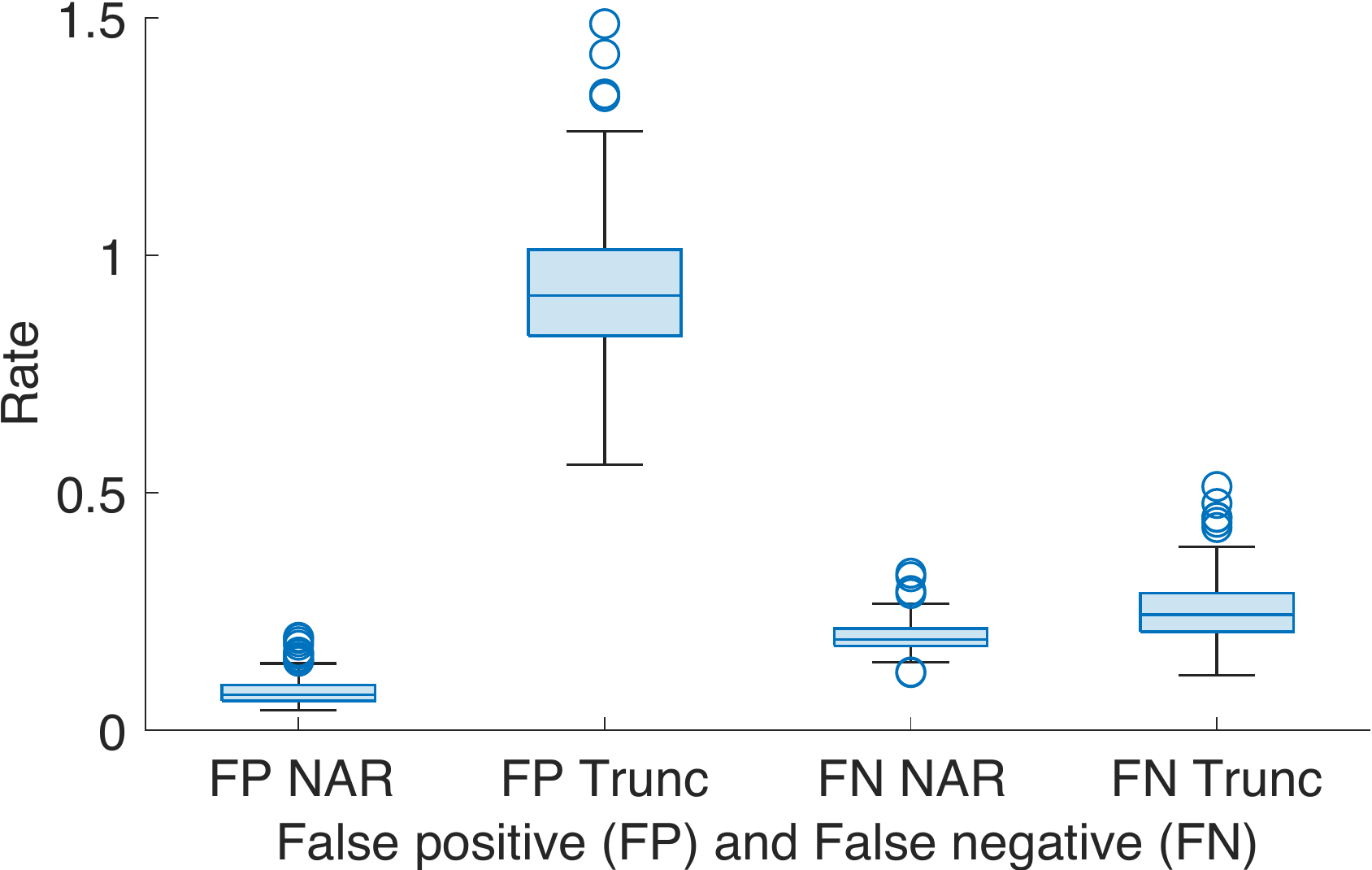}} 
		 \hspace{-2mm}
		  \subfigure[Stochastic force $\sigma=1$]{\includegraphics[width =0.45\textwidth]{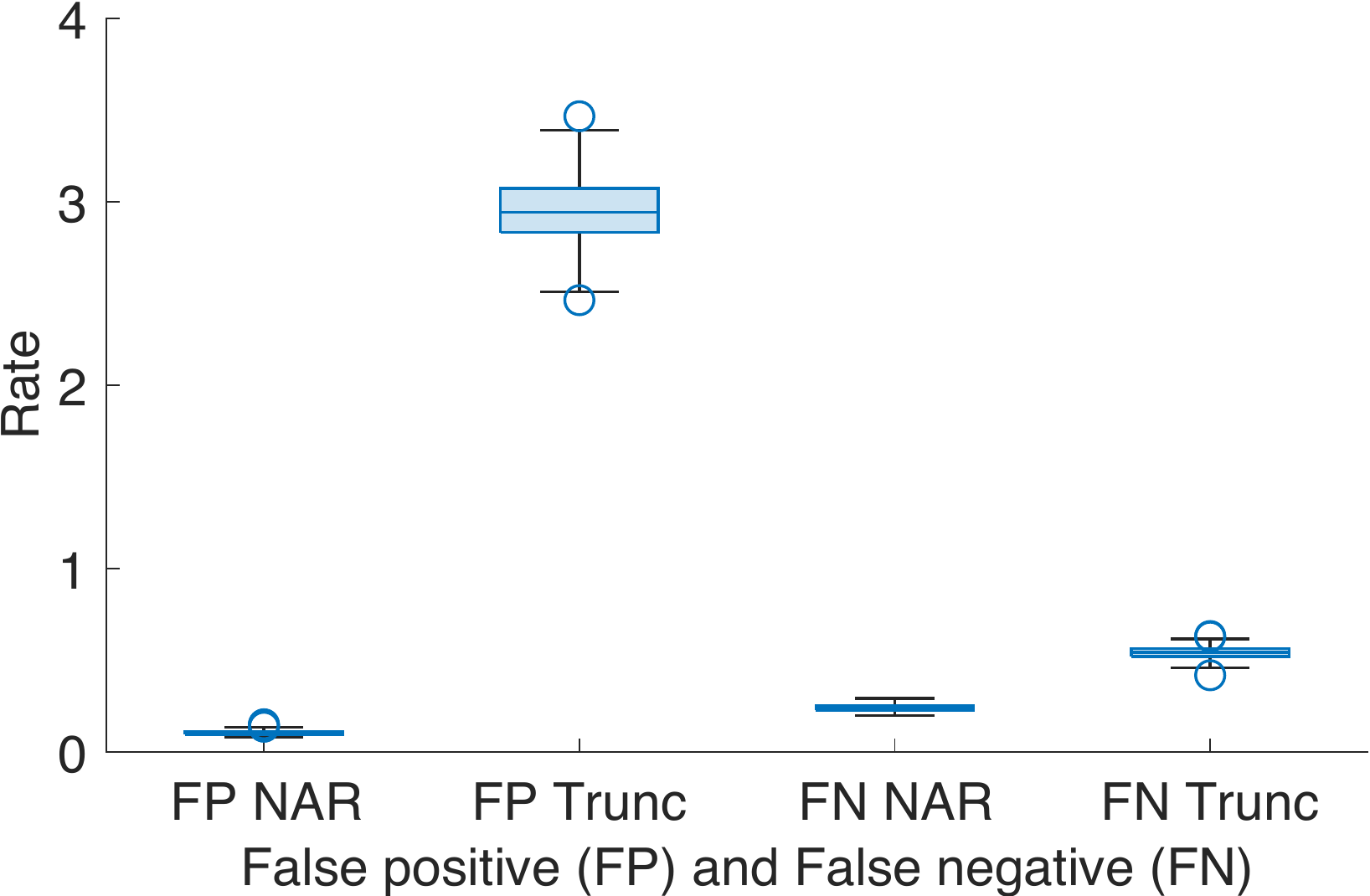} }
		\vspace{-3mm}
\caption{\small Rate of false positive and false negative predictions of shocks in box plots. We compare the NAR model and the truncated model with the $8$-mode projection of the true solution. The rates are out of 200 simulations for each case. In each box, the central mark indicates the median, and the bottom and top edges of the box indicate the 25th and 75th percentiles, respectively. The whiskers extend to the most extreme data points not considered outliers, and the outliers are plotted individually in circle marker symbol.}\label{fig:rateFalse}
\end{figure}

The robust superior performance of the NAR model is further confirmed by the statistics of the false negative and false positive rates. For this purpose, we run the above simulation of the full model and the reduced models for 200 different realizations of the solution of the full model, where the initial conditions are sampled from a long trajectory of the full model. Figure \ref{fig:rateFalse} presents the box plots of the false positive and false negative rates associated with these predictions. The NAR reduced model has significantly smaller false rates than the truncated system for both forcing regimes, and the variation of these rates is also significantly smaller; see the caption of Figure \ref{fig:rateFalse} for more details. 

\subsection{Assimilation and prediction from noisy observations}  \label{Sec_noisy_IC}
Noisy observations are commonly seen in practice, due to either measurement or sampling errors presented in the data collection process. When the observations are noisy, the shock trace prediction consists of two stages: a data assimilation stage to estimate the initial conditions by filtering, and a prediction stage that advances forward in time from the estimated initial conditions. In the data assimilation stage, we filter out the noise in the observations by the ensemble Kalman filter (EnKF) \cite{BLE98,Eve07} (see Appendix \ref{secAppendix:enkf} for a brief review), with respectively the NAR reduced model and the truncated system as the forecast model. 

\paragraph{Test settings.}
For each realization of the true solution, the corresponding shock trace prediction experiment consists of first performing data assimilation using a reduced model (either the NAR model or the truncated system) in the time interval $[0,5]$, and then carrying out prediction using the same reduced model over the time interval $(5,10]$. Thus, we have $500$ time-steps in both stages since the time-step is set to be $\delta=0.01$ for the reduced model.

To generate noisy observations for the data assimilation, we add independent Gaussian noise to the first 8 Fourier modes' real and imaginary parts with a standard deviation $0.01$ and mean zero ( while preserving that $u_{-k} = u_{k}^*$). These noisy perturbations are relatively large, with the ratio between this standard deviation 0.01 and the mean absolute values of the 1st to the 8-th mode of the true dynamics ranging from about 7\% to 40\% for the weak stochastic force regime ($\sigma =0.2$), and from 2\% to 14\% 
for the strong stochastic force regime ($\sigma =1$). This ratio increases as the wave number $k$ increases because the mean absolute value decreases as $k$ increases.

In the assimilation stage, the forecast model of the EnKF (either the NAR model or the truncated system) is initialized from an ensemble of  $100$ initial conditions randomly sampled from a Gaussian distribution centered at the noisy observation at time $t=0$ with standard deviation $0.0254$, which is slightly larger than the standard deviation of the observation noise. Note that the standard deviation used in the construction of these initial conditions is much large than that of the observation noise. In the prediction stage, we simply simulate the forecast model  with the end point from each of the assimilation ensemble as the initial condition. Throughout this section, the NAR model is the same as the one used in Section~\ref{Sec_noise_free_IC}, i.e., the NAR model is trained offline from noiseless data.

\begin{figure}[tbh!]
    \centering
		\subfigure[Stochastic force  $\sigma=0.2$]{\includegraphics[width =0.8\textwidth]{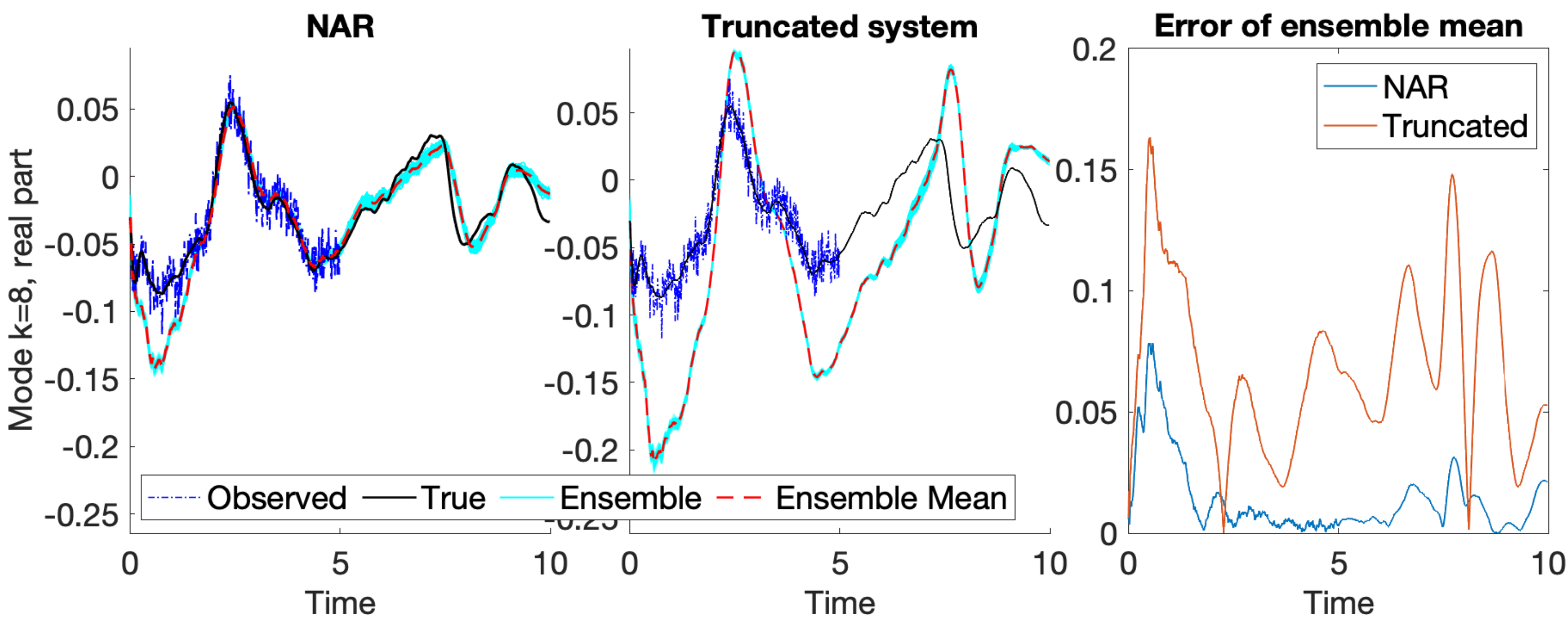}} 
		\subfigure[Stochastic force $\sigma=1$]{\includegraphics[width =0.8\textwidth]{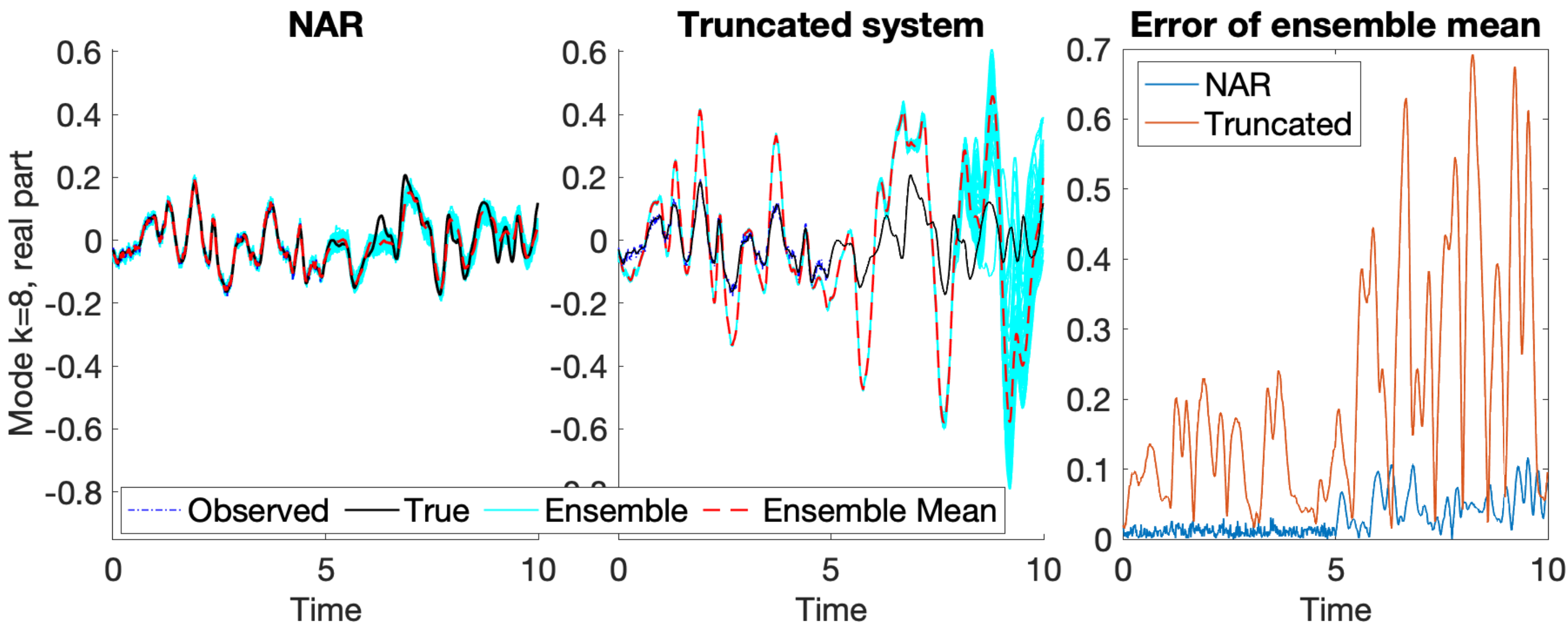} } 
		\vspace{-3mm}
\caption{\small Ensemble trajectories of the real part of the $8$-th Fourier mode in the assimilation stage ($t \in [0, 5]$) and the prediction stage ($t \in (5, 10]$) compared with the corresponding true dyanmics. In both stages and in both cases of weak ($\sigma = 0.2$) and strong ($\sigma = 1$) stochastic forces, the NAR model has ensemble trajectories closer to the true signal (left column)  than the truncated system (middle column), leading to smaller errors in the ensemble mean estimator shown in the right column. 
}\label{fig:stateEst}
\end{figure}

\paragraph{Performance comparison for a typical prediction.} 
As in Section~\ref{Sec_noise_free_IC}, we first compare the performance for a typical realization of the true solution. For this purpose, we show the ensemble trajectories in the data assimilation and the prediction stages. The results are presented in Figure \ref{fig:stateEst} for the real part of the mode with wave number $k=8$, in the settings of a weak stochastic force (top row, $\sigma = 0.2$) and a strong stochastic force (bottom row, $\sigma = 1$).
 
At the data assimilation stage ($t \in [0, 5]$), The EnKF ensembles of the NAR model are much closer to the true trajectory than the truncated system's. Also, the NAR model's ensemble mean estimators have errors less than the observation noise's standard deviation after a short time period, whereas the truncated system's ensembles deviate far away from the truth; see the right column of Figure \ref{fig:stateEst}.
The NAR model's prediction ensembles continue to be spread around the true trajectory, whereas those of the truncated system struggle to make a reasonable prediction for both $\sigma = 0.2$ and $\sigma = 1$.

Similar superiorgood performances are observed for the NAR model in the ensemble prediction of other modes with wave number $k$ less than $8$ at both the data assimilation stage and the prediction stage. The truncated system's prediction skill gradually improves as $k$ decreases, because the signal to noise ratio increases improves as $k$ decreases, but the skill remains inferior to that of the NAR model. Finally, it is worth pointing out that we set the prediction interval to be $(5,10]$ here just for illustration purpose. Since the full model is not chaotic, and the same realization of the stochastic force $f$ is used in both the full model and the reduce systems, one can expect that the good prediction skill of the NAR model shown here to still hold for even longer prediction intervals. 

\begin{figure}[htb]
    \centering 
		\subfigure[Stochastic force  $\sigma=0.2$]{\includegraphics[width =0.75\textwidth]{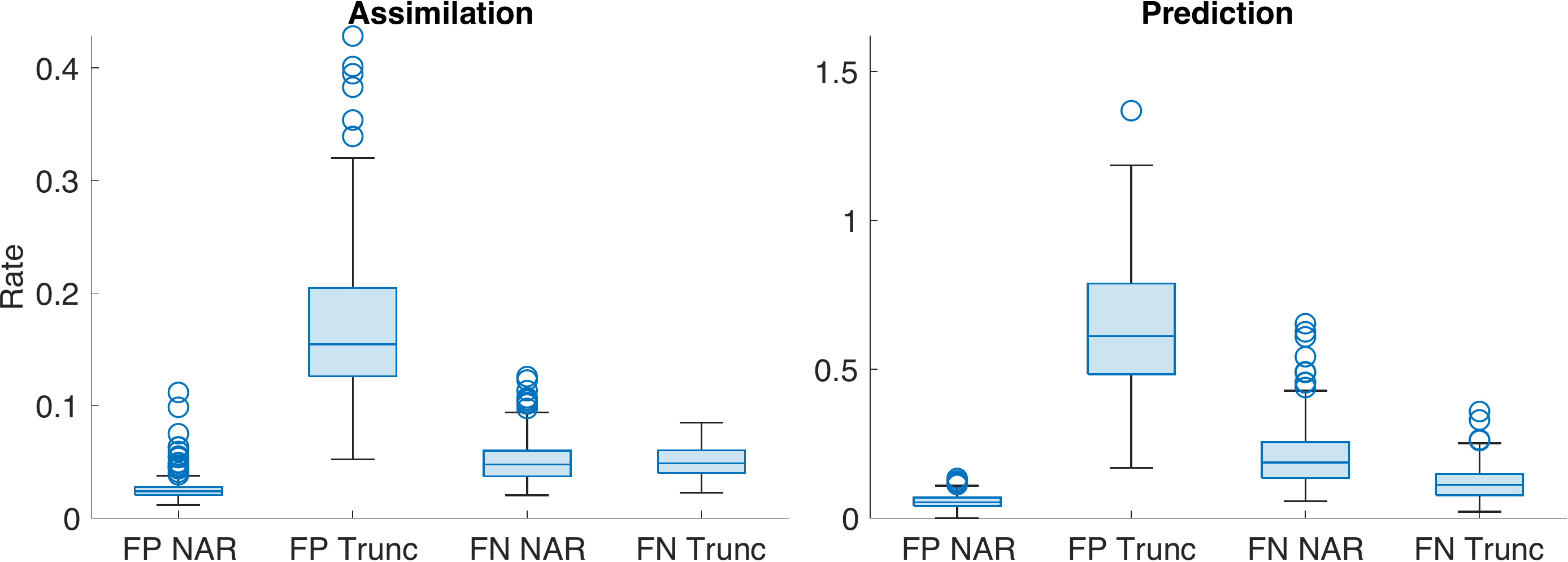}} %
		\subfigure[Stochastic force $\sigma=1$]{\includegraphics[width =0.75\textwidth]{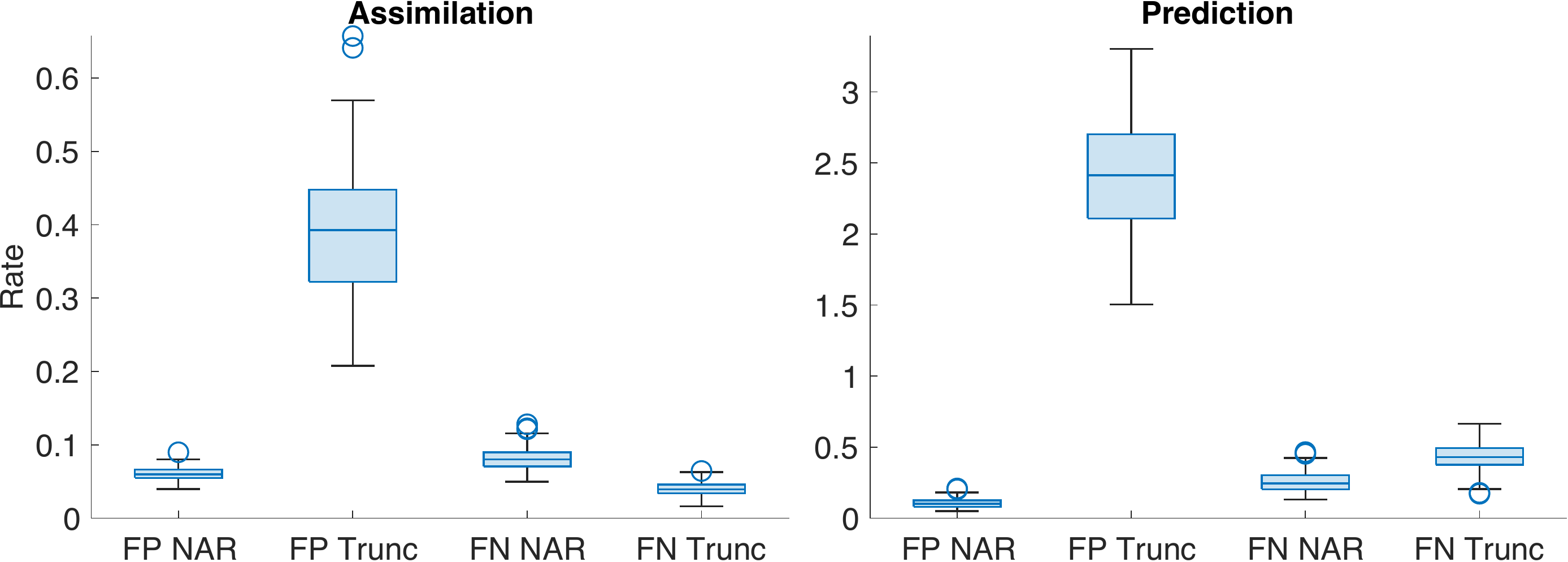} } \hspace{-2mm}
		\vspace{-3mm}
\caption{\small Rate of false positive (FP) and false negative (FN) in the prediction of the binary shock traces obtained from the NAR model and the truncated system.  These rates are computed from the prediction for 200 realizations of the true solution. In each prediction experiment, the spatiotemporal fields for the NAR model and the truncated system used in computing the false rates are taken to be the corresponding ensemble mean of the EnKF with 100 particles. See also the caption of Figure~\ref{fig:rateFalse} for the meaning of each element in the box plots.}\label{fig:rateFalse_DA}
\end{figure}

\paragraph{Performance comparison in multiple simulations.}  To assess the robustness of the NAR model, we repeat the above simulation for 200 different realizations of the true solution and present the performance using again box plots for the rates of false predictions. The procedure of computing the false positive and false negative rates is the same as before, but now carried out separately for the assimilation stage and the prediction stage. For each realization of the true solution, the spatiotemporal fields for the NAR model and the truncated system used in computing these rates are taken to be the corresponding ensemble mean from simulations with 100 different initial conditions, which are sampled in the same way adopted to produce Figure~\ref{fig:stateEst}. 

The NAR model has significantly smaller rates of false predictions than the truncated system in both stages and for both regimes of the stochastic forces. In particular, the improvement is significant in the prediction stage, reducing the median rates of false positive from $0.61$ to $0.05$ in the weak forcing case and from about $2.41$ to less than $0.10$ in the strong forcing case. In the assimilation stage, the improvement is less significant, reducing the median rates of false positive from $0.15$ to  $0.02$ in the weak forcing case and from $0.39$ to $0.06$ in the strong forcing case. The assimilation stage has smaller rates than the prediction stage because of the information supplied by the observation data.

We have thus illustrated that the NAR model is able to predict accurately the pathwise behavior, i.e., the timing and location, of the shocks with uncertainty quantification. This contrasts to predicting only the long-term statistics of the fat-tailed PDFs, which is often the focus of forecasting many other extreme events using reduced models. While the same realizations (after coarsening) of the stochastic force as the full model were applied to the reduced models for diagnostic purpose, the accurate pathwise performance offered by the NAR model presented above shows that even when the stochastic force is unknown in practice, the prediction provided by the $K$-mode NAR model would be comparable with those obtained from the $K$-mode projection of the prediction computed from the full model itself.

\section{Discussion and conclusions} \label{sec:conclusion}
To summarize, this work shows that a reduced model with closure terms, namely a nonlinear autoregression (NAR) model, can accurately predict the timing and locations of random shocks for the viscous stochastic Burgers equation, even though the reduced model cannot represent the precise shocks. A key element is a new characterization of shocks, called \emph{shock trace}, which is a binary indicator that reveals the timing and location of the occurrence of a shock. The shock trace is defined via empirical resolution-adaptive thresholds.

The NAR model describes the evolution of only the leading $K$ Fourier modes with $K$ being a friction of the full model's and with a time-step size ten times the full model's. Thus, it reduces the computational cost orders of magnitude. The NAR model predicts the shock trace almost as good as those computed from the $K$-mode projection of the full model's solution. 
The NAR model consistently outperforms the corresponding Galerkin truncated system in forecasting the shock trace from both noiseless and noisy observations. Thus, the data-driven closure in the NAR model plays an essential role in shock trace prediction.

We conclude by reiterating some key elements that contributed to the success of the closure model utilized here, and also mention potential future works.  

\medskip
\noindent{\bf Data representation and importance of closure terms.} 
The reduction and prediction approach presented here is independent of the orthogonal basis adopted to decompose the solution field. Instead of the Fourier basis, one could use for instance empirically computed modes obtained from {\it e.g.}~the proper orthogonal decomposition method~\cite{HLB96} and its variants \cite{taira2017modal} or the dynamic mode decomposition \cite{rowley2009spectral,schmid2010dynamic}. 

The superior performance of the NAR model over the truncated system demonstrates the importance of incorporating suitable closure terms into the reduced model in order to effectively represent the nonlinear feedback from the unresolved-scale variables to the resolved ones. Such nonlinear feedback effects are important to be properly approximated here since a significant amount of energy in the shocks can be contained in the high-frequency variables which are not resolved by a low-dimensional reduced system (cf.~Table~\ref{tab:energy_frac}). Of course, besides the NAR model, many other closure modeling techniques can potentially be used for this purpose; see e.g.~\cite{majda2001mathematical,LinLu21,wouters2013multi,xie2018data,ahmed2021closures,chekroun2020variational}. 

\medskip
\noindent{\bf Predicting key information of extreme events by reduced models.} By relaxing the goal from predicting the full shock profiles to predicting only their space-time locations, we arrive at a task within the representation capacity of a low-dimensional reduced subspace, which is thus much more accessible by a reduced-order modeling approach.  The indicator function \eqref{eq:binaryShockTrace} introduced here for identifying binary shock traces is based on a simple resolution-adaptive threshold defined by \eqref{eq:tau}, which is computed empirically using an ensemble of the true solution trajectories or their projections onto a given number of Fourier modes. These binary shock traces encode the shock locations in space-time, providing thus important partial information about the shocks. 

While we illustrated the approach on the stochastic Burgers equation, it is expected that the developed strategy will be of interest for data-driven predictive modeling of a broad range of extreme events, particularly those beyond the representation capacity of a reduced order modeling framework. Note that the concerned extreme events do not have to be the extreme values \cite{ghil2011extreme,kwasniok2012data,qi2020using,wan2018data} of the states of the model; the shocks correspond to extreme values in the spatial gradient field instead. It is reasonable to expect that the representation capacity offered by a fixed subspace will be violated more often when higher derivatives are involved in defining the concerned extreme events.

\medskip
\noindent{\bf Prediction of the exact shocks}. Going beyond shock trace prediction, if one would wish for instance to approximate the sharp gradients presented in the shock profiles, it would require either to resolve more Fourier modes, or switch to a data-adaptive empirical basis such as those mentioned above, or adopt a combination of good low-dimensional reduced models with additional techniques to recover the unresolved high-frequency modes. A few possibilities are available for the latter option depending on the setup. 

For instance, when observations are only available for the low-frequency modes, one can design computationally efficient data assimilation strategies within the conditional Gaussian framework \cite{chen2018conditional,chen2021conditional} to approximate the dynamics of the high-frequency modes with quantified uncertainties by a suitable dynamical model for the unresolved modes. On the other hand, one can parameterize the high frequency modes using suitable random functions of the resolved low frequency dynamics. These random functions can be designed for instance through the dynamics-based parameterizing manifold approach \cite{CLW15_vol2,chekroun2020variational} or other data-driven/machine learning approaches. When good parameterizations of the unresolved modes are available, they can be used to both build reduced models for the resolved modes and provide approximations of the unresolved modes.

\section*{Acknowledgments}
The research of N.C. is partially funded by the Office of VCRGE at UW-Madison and ONR N00014-21-1-2904. The work of H.L. is partially funded by NSF Award DMS-2108856. The work of F.L. is funded by NSF Award DMS-1821211 and DMS-1913243.

\appendix
\section{The ensemble Kalman filter} \label{secAppendix:enkf}

The Ensemble Kalman Filter (EnKF) is a Monte Carlo implementation of
Bayesian filtering with the Kalman filter update \cite{Eve94,EVL96,HM98,BLE98}. It uses an
ensemble of random samples, also called particles, to approximate the forecast and analysis distributions by Gaussian distributions whose means and covariances are given by ensemble means and covariances. Among various EnKF algorithms, we consider the version with perturbed observations, introduced in \cite{BLE98,HM98}, and we refer to \cite{LBS10} for a comparison of different versions of EnKF algorithms.

Suppose the filter uses a forecast model
\begin{equation} \label{frmodel}
\mathbf{x}_{n}=\mathbf{F}_{n}(\mathbf{x}_{n-l:n-1}),
\end{equation}
where $\mathbf{x}_{n}\in \mathbb{R}^{d_{x}}$ is the state variable,  $\mathbf{x}_{n-l:n-1} = (\mathbf{x}_{n-l}, \dots, \mathbf{x}_{n-1})$, and $\mathbf{F}_{n}$ is a
forecast operator at time $n$ which maps $\mathbb{R}^{l\times d_{x}}$ to $\mathbb{R}^{d_{x}}$ with $1\leq l \leq n-1$. The forecast model can be either
stochastic or deterministic, and either Markovian (e.g. $l=1$) or non-Markovian (e.g. $l>1$). The state variable is observed through a linear observation operator with Gaussian noise:
\begin{equation*}
\mathbf{z}_{n}=H\mathbf{x}_{n}+\mathbf{\epsilon}_{n},
\end{equation*}%
where $H\in \mathbb{R}^{d_{z}\times d_{x}}$ is the observation matrix, and the $%
\mathbf{\epsilon}_{n}\sim N(0,R)$ are independent Gaussian noises. {\color{black}In this study, we assume that the observation matrix $R$ is known. }

The EnKF iterates the following two steps, with an initial ensemble of
particles $\{\mathbf{x}_{0}^{a,(i)}, i=1,\dots ,M\}$ sampled from the forecast distribution of the state variable $\mathbf{x}$ (e.g. the stationary distribution of the forecast model).
\begin{enumerate}
\item Forecast step: from the ensemble $%
\{\mathbf{x}_{1:n-1}^{a,(i)}\}$ at time $n-1$, generate a forecast ensemble $\{\mathbf{x}
_{n}^{f,(i)}\}$ using the forecast model in (\ref{frmodel}), i.e. $\mathbf{x}_{n}^{f,(i)}=F_{n}(\mathbf{x}_{n-l:n-1}^{a,(i)})$. Here the superscript in $\mathbf{x}_{n}^{f}$ denotes the ensemble from the forecast model, and the superscript in $\mathbf{x}_{n}^{a}$ denotes the ensemble of the posterior distribution after assimilating data in the following analysis step.
If the forecast model is stochastic, independent realizations should be used at different times.

\item Analysis step. Given new observation $\mathbf{z}_{n}$,
update the forecast ensemble to get a posterior ensemble of $\mathbf{x}_{n}$,
\begin{equation} \label{EnKF_update}
\mathbf{x}_{n}^{a,(i)}=\mathbf{x}_{n}^{f, (i)}+K_{n}(\mathbf{z}_{n}^{(i)}-H\mathbf{x}_{n}^{f,(i)}),
\end{equation}%
for $i=1,\dots ,M$, where the Kalman gain matrix is
\begin{equation} \label{EnKF_K}
K_{n}=C_{n}^{f}H^{T}(HC_{n}^{f}H^{T}+R)^{-1},
\end{equation}%
where the matrix $C_{n}^{f}$ is the sample covariance of the forecast ensemble:%
\begin{equation*}
C_{n}^{f}=\frac{1}{M-1}\sum_{i=1}^{M}\left( \mathbf{x}_{n}^{f,(i)}-\overline{\mathbf
{x}}_{n}^f\right) \left( \mathbf{x}_{n}^{f,(i)}-\overline{\mathbf{x}}_{n}^f\right) ^{T},
\end{equation*}%
where $\overline{\mathbf{x}}_{n}^f=\frac{1}{M}\sum_{i=1}^{M}\mathbf{x}_{n}^{f,(i)}$ and the $\mathbf{z}_{n}^{(i)}$ are obtained by adding random perturbations $\mathbf{\epsilon}_{n}^{(i)}\sim N(0,R)$ to $\mathbf{z}_{n}$:%
\begin{equation*}
\mathbf{z}_{n}^{(i)}=\mathbf{z}_{n}+\mathbf{\epsilon}_{n}^{(i)}.
\end{equation*}
\end{enumerate}

\section{Parameters in the NAR models}
The parameters for the NAR models in the regimes of weak and strong stochastic forcing are shown in Table \ref{tab:parNAR_sigma_p2} and \ref{tab:parNAR_sigma1} (see Section~\ref{Sec_model_selection} for the specifications). 

\begin{table}[H]\vspace{-3mm}  
\centering
\caption{Parameters in NAR model: weak stochastic force $\sigma=0.2$}\label{tab:parNAR_sigma_p2}
\begin{tabular}{c | cccc cccc }
 k & 1 & 2 & 3& 4 & 5 &6 &7 & 8 \\
\hline
$c^v$  & -0.01 & -0.05 & -0.12 & -0.26 & -0.73 & -1.28 & -2.02 & -2.68  \\
$c^R$ & 1.08 & 1.04 & 1.00 & 0.93 & 0.93 & 0.83 & 0.65 & 0.26  \\
$c^f$ ($\times 10^{-3}$)  & -0.12 & -0.47 & -1.07 & -1.90 & -0.00 & -0.00 & -0.00 & -0.00  \\
$c^w$ ($\times 10^{-5}$)& -1.38 & 4.59 & -3.41 & -4.47 & -5.48 & -5.37 & -0.00 & -5.91  \\
$\sigma_g$  &0.04 &0.13 &0.23 &0.32 &0.44 &0.55 &0.70 &0.93
   \end{tabular}
   \end{table}
\vspace{-3mm}   
   \begin{table}[H] \vspace{-3mm}  
\centering
\caption{Parameters in NAR model: strong stochastic force $\sigma=1$}\label{tab:parNAR_sigma1}
\begin{tabular}{c | cccc cccc }
 k & 1 & 2 & 3& 4 & 5 &6 &7 & 8 \\
\hline
$c^v$ & -0.05 & -0.23 & -0.54 & -1.07 & -2.71 & -4.49 & -6.84 & -9.20 \\ 
$c^R$ & 1.08 & 1.02 & 0.97 & 0.87 & 0.88 & 0.78 & 0.64 & 0.35 \\
$c^f$ ($\times 10^{-3}$) & -0.26 & -1.05 & -2.39 & -4.25 & 0.00 & -0.00 & 0.00 & -0.00  \\
$c^w$ ($\times 10^{-6}$)  & -0.70 & -6.88 & -1.87 & -4.57 & -1.84 & -2.27 & -1.43 & -7.74 \\ 
$\sigma_g$  &0.92 &2.06 &3.19 &4.08 &5.12 &5.96 &7.28 &9.74 
   \end{tabular}
   \end{table}

\bibliographystyle{myplain}
\bibliography{ref_FeiLU2021_11,Model_reduction,ref_Burgers,ref_chen}

\end{document}